\documentclass[11pt]{amsart}  % use ``amsart'' instead of ``article'' for AMSLaTeX format
\usepackage[left=1.2in,right=1.2in,bottom=1.2in]{geometry}  
\usepackage[T1]{fontenc}

\usepackage{amssymb,amsmath,amscd}
\usepackage{amsthm}
\usepackage[parfill]{parskip}
\usepackage{upquote}
\usepackage{upgreek}
\usepackage{mathrsfs}
\usepackage{latexsym}
\usepackage{graphicx} % Required for inserting images
\usepackage[x11names]{xcolor}
\usepackage{tikz-cd}
\usepackage{tikz} \usetikzlibrary{positioning,arrows.meta}

\usepackage{enumitem}
\setlist{itemsep=0.2em, parsep=0.2em} % Adjust these values as needed
\setlist[itemize,1]{label=\ensuremath{\blacktriangleright}}
\setlist[itemize,2]{label=\ensuremath{\triangleright}}

\theoremstyle{plain}
\newtheorem{theorem}{Theorem}[section]
\newtheorem{lemma}[theorem]{Lemma}
\newtheorem{proposition}[theorem]{Proposition}

\theoremstyle{definition}
\newtheorem{definition}[theorem]{Definition}
\newtheorem{example}[theorem]{Example}
\newtheorem{construction}[theorem]{Construction}

\theoremstyle{remark}
\newtheorem{remark}[theorem]{Remark}

\DeclareMathOperator{\coker}{coker}                 % Cokernel of a map

\newcommand{\torsor}{\mathcal{P}} % TORSOR
\newcommand{\glosection}{\sigma} % GLOBAL SECTION
\newcommand{\res}{\mathbf{r}} % RESTRICTION
\newcommand{\Aut}{\mathrm{Aut}} % AUTOMORPHISM GROUP
\newcommand{\face}{\to} % FACE OPERATION, FORMERLY \lhd
\newcommand{\vertexneg}{u} % LEFT VERTEX
\newcommand{\vertexpos}{v} % RIGHT VERTEX

\newcommand{\vect}[1]{\mathbf{#1}}
\title[Torsors \& Visual Paradox]{Obstructions to Reality : \\
Torsors \& Visual Paradox}
\author{Robert Ghrist}
\address{Departments of Mathematics and Electrical \& Systems Engineering \\ University of Pennsylvania \\ Philadelphia, PA 19104}
\email{ghrist@math.upenn.edu}
\author{Zoe Cooperband}
\address{Department of Electrical \& Systems Engineering \\ University of Pennsylvania \\ Philadelphia, PA 19104}
\email{zcooperband@gmail.com}

\begin{document}

\begin{abstract}
Visual paradoxes like the Penrose staircase present a fundamental tension: locally coherent geometric relationships that cannot be realized globally. Inspired by Penrose's observations connecting such paradoxes to cohomology, we develop a mathematical framework that characterizes this phenomenon through network torsors and sheaf cohomology. Network torsors capture the essential nature of visual paradoxes by formalizing relative geometric attributes (height changes, orientation flips) without requiring absolute measures. We demonstrate that a significant class of visual paradoxes can be characterized as non-trivial network torsors, with their obstruction to global consistency quantified by elements of first cohomology $H^1$. This framework enables analysis of classical paradoxes and construction of novel examples on various topological spaces. Key contributions include the first visual paradox with nonabelian holonomy (the Klein ladder, whose holonomy takes values in the infinite dihedral group), and paradoxes driven by boundary conditions rather than loops, analyzable via non-constant structure sheaves and relative cohomology. Our approach unifies diverse visual paradoxes under a single mathematical principle: the obstruction to globalizing locally consistent geometric relationships.
\end{abstract}


\maketitle

%%%%%%%%%%%%%%%%%%%%%%%%%%%%%%%%%%%%%%%%%%%%%
\section{Introduction}
\label{sec:intro}
%%%%%%%%%%%%%%%%%%%%%%%%%%%%%%%%%%%%%%%%%%%%%

\begin{figure}[hb]
    \centering
    \setlength{\fboxsep}{0pt} % optional: no padding between box and image
    \fbox{\includegraphics[height=6.5cm]{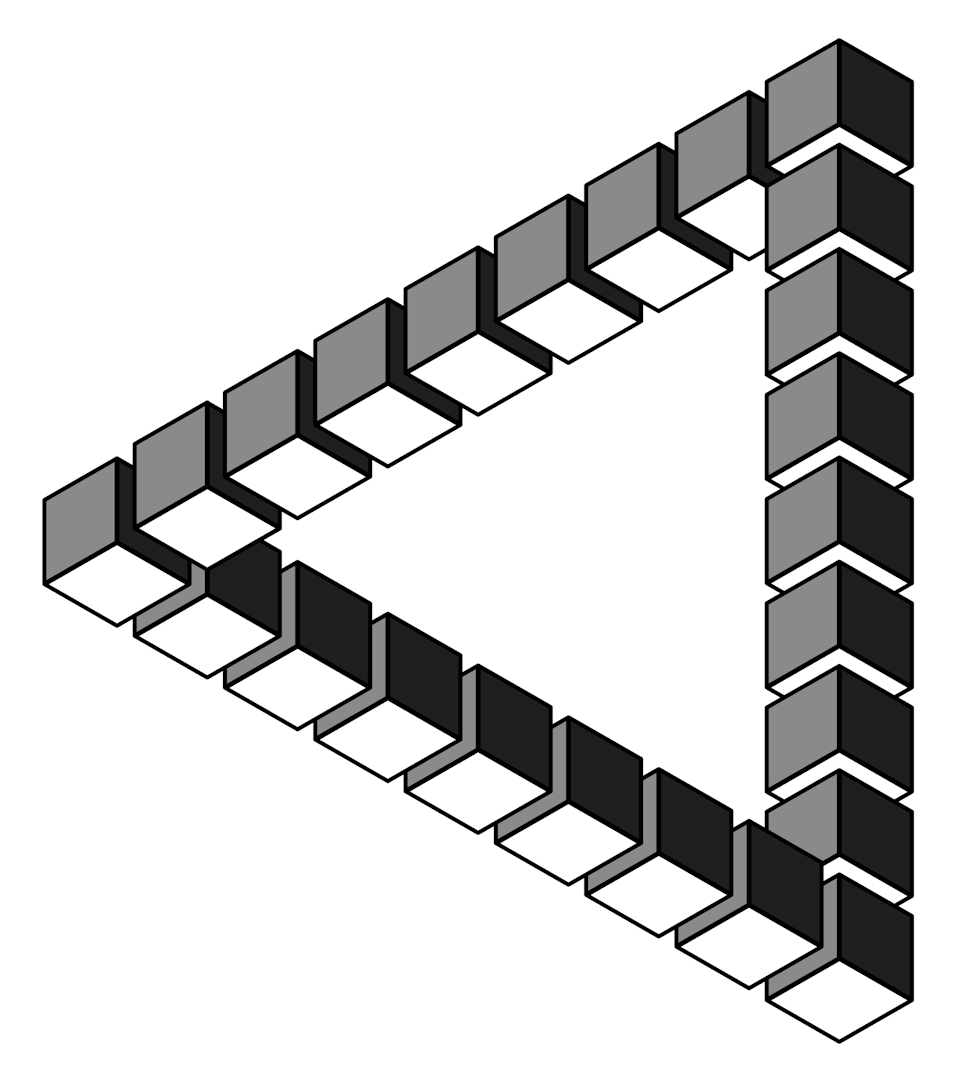}
}%
    \hspace{0.5em}%
    \fbox{\includegraphics[height=6.5cm]{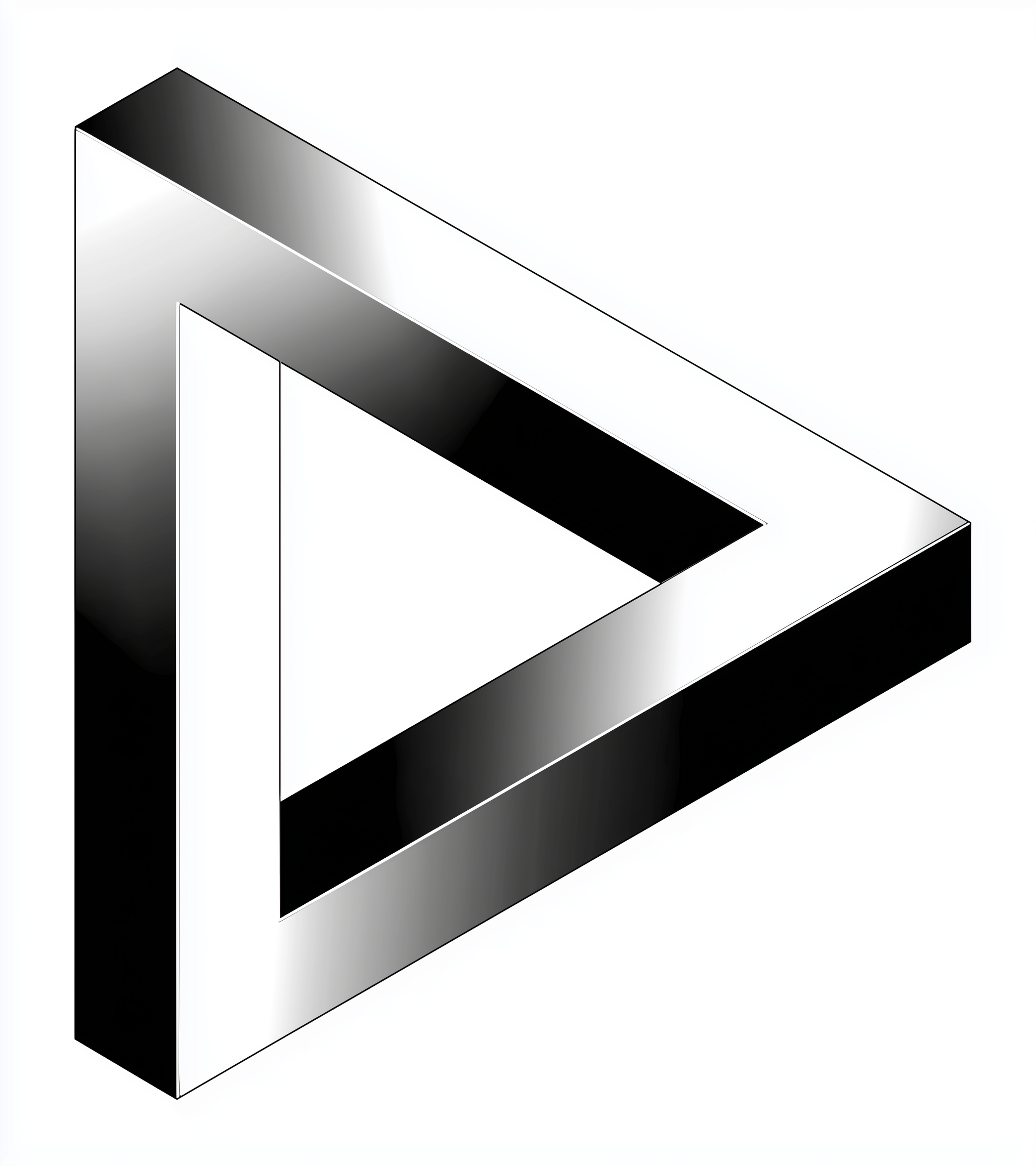}}    \caption{\small [left] The impossible triangle of Reutersv\"ard (1934); [right] the classic Penrose triangle (1954--56). Each is locally coherent, yet no 3-D structure realizes the local readings of all its corners and beams at once: forced-perspective sculptures realize the image, not those readings.}
    \label{fig:classics}
\end{figure}

\begin{figure}
    \centering
    \setlength{\fboxsep}{0pt}
    \fbox{\includegraphics[height=5cm]{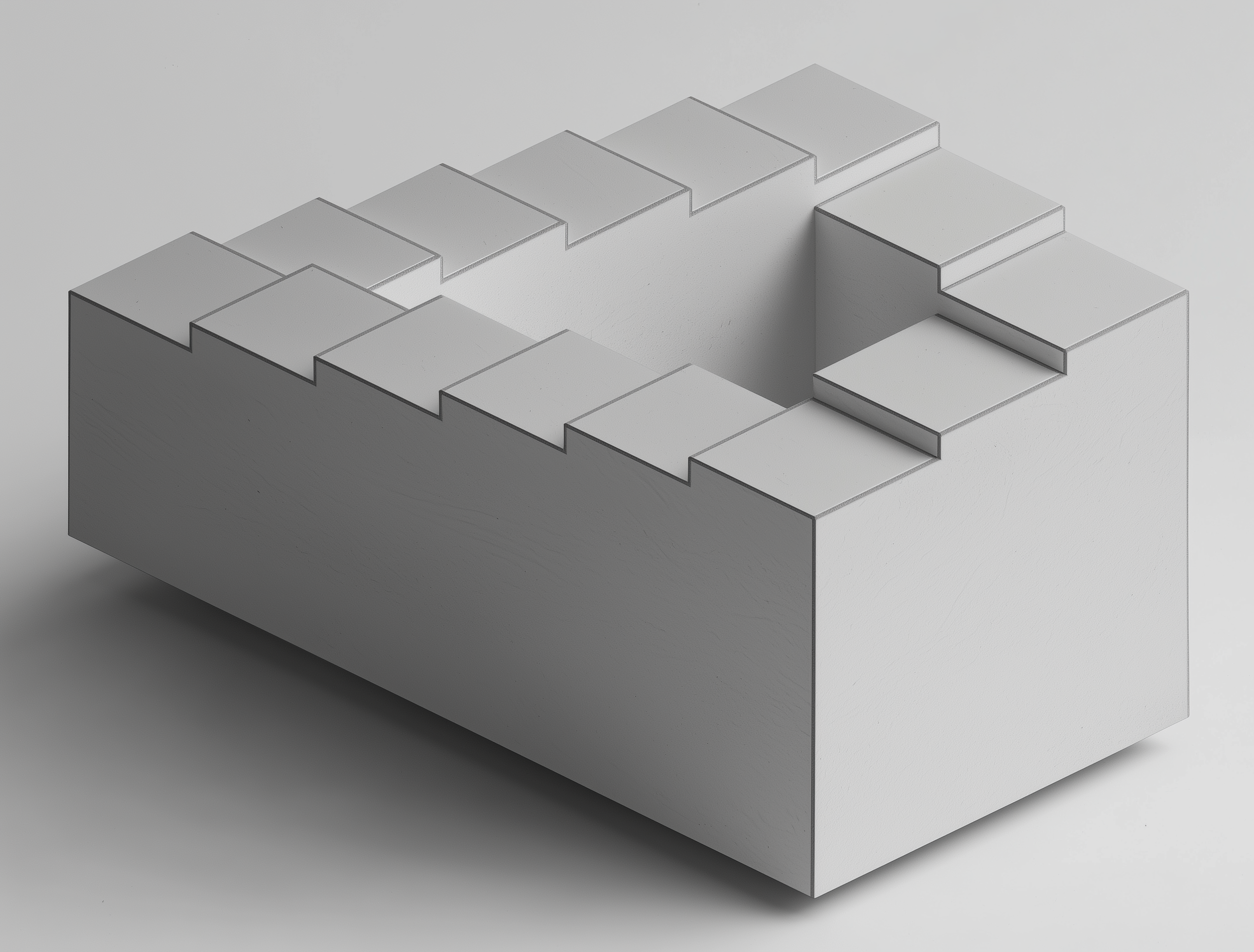}}
    \caption{\small The Penrose staircase (1954-56) has paradoxical height as one traverses the loop.}
    \label{fig:classics2}
\end{figure}

Impossible figures represent a fascinating intersection of visual perception, geometry, and topology. Epitomized by the Penrose triangle or staircase and elaborated in the works of Reutersv\"ard, Escher, and others, these structures present locally coherent geometric information -- segments of stairs appear correctly joined, perspectives seem locally valid -- yet defy global realization within standard Euclidean space. The observer perceives a paradox: a structure that appears logically consistent piece by piece, yet fundamentally impossible when considered as a whole. Figure \ref{fig:classics} recalls the Penrose triangle (or {\em tribar}) and Figure \ref{fig:classics2} gives an example of a Penrose staircase \cite{penrose1958impossible}.

\begin{figure}[htb]
    \centering
    \setlength{\fboxsep}{0pt}
    \fbox{\includegraphics[height=5cm]{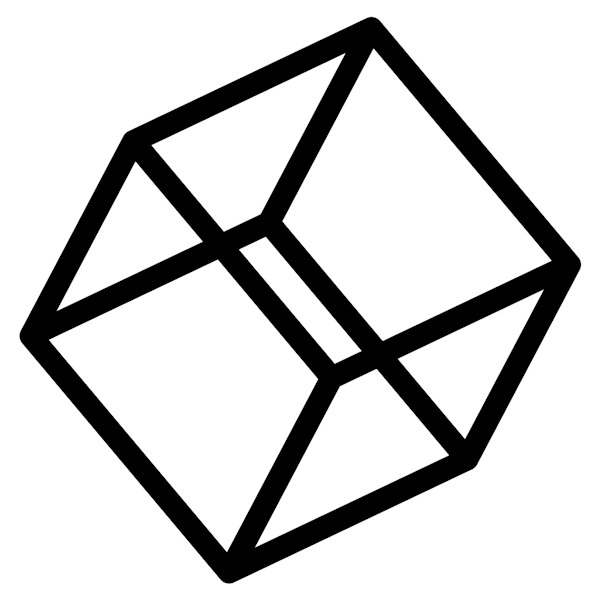}}
    \hspace{0.5em}
    \fbox{\includegraphics[height=5cm]{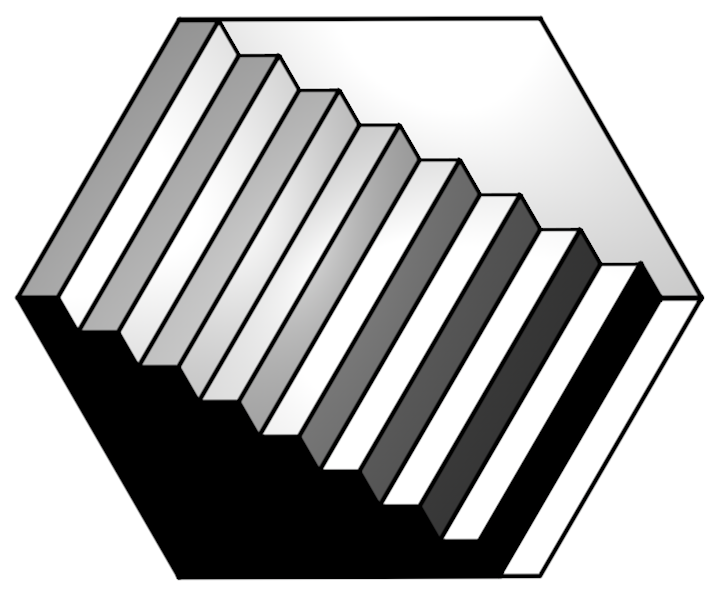}}
    \caption{\small The Necker cube (1832) [left] and the Schr\"oder staircase (1858) [right] each have a bistable {\em gestalt switch} that can serve as the basis for impossible figures, see, {\em e.g.,} Figures \ref{fig:necker_gradient} and \ref{fig:boundary}.}
    \label{fig:classics3}
\end{figure}

The core of the paradox lies in a conflict between local geometric cues and global topological constraints. In the Penrose staircase, each segment appears to ascend or descend in a consistent direction, yet traversing the complete circuit returns one to the starting point without a net change in elevation. Another class of impossible figures -- commonly seen in the works of Escher and other popular artists -- relies on a {\em gestalt switch} in which an ambiguous image can be held (with effort) in tension between two stable perceptions. Figure \ref{fig:classics3} gives the classic examples of the Necker cube and the Schr\"oder staircase. By manipulating these switches with some global constraints, one can generate paradoxical figures without a consistent global interpretation.

The connection between impossible figures and cohomology was recognized in Penrose's 1992 paper \cite{Penrose1992Cohomology}, which identified the relation to the first cohomology group of the circle $H^1(S^1 ; \mathbb{R}^+)$. He also included the Necker cube as an example with a paradox detected by sheaf cohomology with coefficients in $\mathbb{Z}_2$. While the paper is light on details, the insight that cohomology is a natural framework for analyzing the global inconsistencies that emerge from locally coherent visual cues is seminal. His paper ends with the following intriguing lines.

\small
\begin{quote}
{\em ``More complicated figures with "multiple impossibilities" (see, e.g.
[Penrose \& Penrose, 1958]) can also be analyzed in this way, but for this we should require a more complete description of what a (Cech) cohomology group actually is... I believe
that considerations such as these may open up intriguing possibilities for further exotic types of impossible figure. I hope to be able to consider such matters at a later date.''}
\\
--- R. Penrose \cite{Penrose1992Cohomology}
\end{quote}
\normalsize

This paper aims to return to Penrose's goal of generating more exotic types of visual paradoxes, and with it develop a rigorous mathematical framework for understanding and analyzing such figures. Our central thesis is that a significant class of geometric visual paradoxes can be precisely characterized as non-trivial network torsors: once a model specifies the local states of a figure and how they are compared, the obstruction to a globally consistent interpretation is measured by sheaf cohomology.

%------------------------------------------
\subsection{Towards Torsors}
\label{sec:abs-vs-rel}
%------------------------------------------

The toolkit we propose is motivated by the distinction between two types of visual attributes:

\emph{Absolute attributes} are properties that can be consistently assigned at each point of a figure without reference to other points. Coordinates, color, texture, or material properties exemplify such attributes. When viewing an impossible figure, these absolute attributes cause no cognitive dissonance -- we can consistently see the staircase steps as having certain colors, for example, without any contradiction.

In contrast, \emph{relative attributes} are fundamentally relational, defined only in terms of their alignment to neighboring points. Relative height, orientation, or depth belong to this category. What makes these attributes special is their dependent nature -- they describe how one part of a figure relates to another.

Whether an attribute is absolute or relative is a feature of the model, not of the attribute itself. The Necker state is relative in the middle of the sequence of Figure \ref{fig:necker_gradient}, where only comparisons between neighbors are given, and absolute at its shaded ends, where visual cues pin it.

Visual paradox often emerges precisely when we attempt to globalize these local relations. In the Penrose staircase, each step stands in a well-defined height relationship to its neighbors, yet traversing the complete circuit returns one to the starting point without a net change in elevation -- an impossibility.

Relative attributes are described mathematically by local data-transition functions between adjacent parts of the figure formalized as cocycles. A visual paradox emerges when these cocycles cannot be integrated into a globally consistent assignment for the attribute in question ({\em e.g.}, absolute height, absolute orientation). The obstruction to finding such a consistent global assignment, {\em i.e.}, a global section of the relevant mathematical structure, is precisely what characterizes the paradox. This obstruction is measured by a non-trivial cohomology class. Absolute attributes, like color, typically correspond to data for which a global section is straightforward and do not lead to paradox.

Network torsors provide the ideal mathematical structure to capture this tension. Informally, a torsor can be understood as a collection of relative transformations between points without a canonical choice of absolute positioning. A torsor may still admit global sections -- a trivial torsor admits many, none of them preferred. When a torsor is non-trivial, the obstruction to finding a global section expresses the impossibility of assigning absolute values that respect the relative relationships. This perspective not only clarifies why certain visual attributes create paradoxes while others do not, but also provides a rigorous and algebraic description.

Throughout, the models live on graphs, with prescribed transport along edges and coefficient groups treated algebraically. Each model specifies its base, structure group, local states, and transport, and a claim of impossibility is a claim about that model. When a figure lives on a surface, data on its constraint graph must be shown to extend across the 2-cells; obstructions carried by the 2-cells themselves, in $H^2$, are considered in \cite{GhristGhrist2026Degrees}.

%------------------------------------------
\subsection{Contributions and Organization}
\label{sec:contributions}
%------------------------------------------

\begin{enumerate}
\item We introduce the concept of a network torsor (Definition \ref{def:torsor}), providing a discrete analogue of the classical torsor notion that is suited to the analysis of visual paradoxes, together with interpretation problems in which the local states need not form a torsor (Definition \ref{def:interpretation}). We then interpret a significant class of visual paradoxes as non-trivial network $G$-torsors.

\item We analyze novel paradoxes constructed on surfaces beyond the circle, including the M\"obius strip, projective plane, torus, and Klein bottle. In so doing, we present what is, to our knowledge, the first example of a paradoxical figure with nonabelian holonomy: a {\em Klein ladder}, built with an infinite dihedral torsor and illustrated as an identification space, in which the two orders of traversal of its generating loops end at different heights.

\item We introduce paradoxes driven by boundary constraints without an underlying ``loop'' in the base space. This is handled via a torsor with a non-constant {\em structure sheaf} and relative cohomology. This extends the framework to figures that may have simple underlying topology but exhibit paradox due to conflicting boundary conditions.

\item We describe how pullback and change of coefficients relate models of the same figure or of visually distinct figures, and we return to the Penrose triangle through affine holonomy, which places the height paradoxes of this paper in cohomology with twisted coefficients.
\end{enumerate}

%------------------------------------------
\subsection{Related Work and Context}
\label{sec:work}
%------------------------------------------

There have been numerous attempts to roughly classify or partition different impossible figures using a variety of mathematical notions \cite{Harris1973, Cowan1974, Cowan1977, Kulpa1983}. Cowan's use of braids is close in spirit to the present work: corners and sides of polygonal Penrose-type figures are encoded as braids, and closing products of braids generate possible and impossible figures \cite{Cowan1974}. Each of these figures is a single loop, so any holonomy it carries generates a cyclic group. Cohomological investigation of impossible figures begins notably with Penrose's work on the impossible triangle and Necker cube \cite{Penrose1992Cohomology}, in which he defines an {\em ambiguity group} which encodes the paradox algebraically and a coefficient group. 

Subsequent developments have diverged. One significant branching is in the mathematical study of contextuality, particularly in quantum mechanics. Inspired by Penrose's impossible figures, Abramsky and Brandenburger \cite{AbramskyBrandenburger2011Sheaf} developed a comprehensive sheaf-theoretic framework for analyzing contextual phenomena, where measurements or observations depend on their context. While their focus was primarily on quantum mechanics, they established the foundation for interpreting various paradoxes through the lens of sheaf theory. In particular, the absence of a global section of a presheaf over contextual measurements precisely mirrors the impossibility of global consistency in paradoxical figures.

This connection was further developed by Abramsky et al.~\cite{AbramskyBarbosa2015Paradox}, who linked the cohomological perspective on quantum contextuality with logical paradoxes. Cervantes and Dzhafarov \cite{CervantesDzhafarov2020Impossible} applied contextuality analysis to impossible figures, treated as systems of epistemic random variables, and measured their degree of contextuality. Car\`u \cite{Caru2017Cohomology} generalized the cohomological obstructions of contextuality to higher cohomology groups, organized the resulting notions in a hierarchy, and described first cohomology in terms of torsors for an abelian presheaf; his examples are contextual scenarios, not visual figures.

Few of these works connect explicitly to the visual or perceptual paradoxes of Penrose's initial vision; those that do \cite{CervantesDzhafarov2020Impossible} revisit classical paradoxical images without introducing novel paradoxical figures. There is an important unpublished talk by Sheehy which explored the cohomology of impossible figures, reinterpreting Penrose's work from a computational viewpoint of convex lifts and folds.\footnote{Recorded and available at {\tt https://www.youtube.com/watch?v=CQ7zTx3xJZM}} Sheehy connected the paradox to obstructions in lifting the 2-D figure to 3-D, arguing that cycles in local geometric or occlusion orderings correspond to the non-trivial $H^1$ obstructions. This perspective is extended in the thesis of Cooperband \cite{cooperband2024cellular} where such cycles are explicitly related to liftings via a long exact (co)homology sequence.

A separate line of work treats the realizability of line drawings computationally, by labeling schemes \cite{Huffman1971} and by linear algebra \cite{Sugihara1986}; the kernel-of-projection argument of Section \ref{sec:cubic_staircases} is a small instance of the latter viewpoint. Our work builds upon these foundations while making the contributions outlined in the previous subsection. For the reader unfamiliar with sheaves and torsors, the gallery below is a gentle motivation for the mathematics to follow.

%%%%%%%%%%%%%%%%%%%%%%%%%%%%%%%%%%%%%%%%%%%%%
\section{A Gallery of Paradoxical Stairs}
\label{sec:gallery}     
%%%%%%%%%%%%%%%%%%%%%%%%%%%%%%%%%%%%%%%%%%%%%

Before developing formal mathematical machinery, we present a collection of paradigmatic impossible figures, establishing the visual phenomena we seek to explain and classify. Rather than deal with all possible illusions, we restrict attention to stair-like diagrams and images with forced perspectives.\footnote{One finds in popular artworks any number of paradoxical images depending on orientation, hidden lines, and other mechanisms. The analysis in this paper is not unrelated to such, but, for simplicity, we do not strive for universal applicability.} 

Many of these ``impossible figures'' are novel and based on identification spaces, where lifting to a universal cover can be helpful in visualization.

%------------------------------------------
\subsection{Reutersv\"ard Stairs}
\label{sec:gallery-stairs}
%------------------------------------------

\begin{figure}
    \centering
    \begin{minipage}[t]{0.24\textwidth}
        \centering
        \includegraphics[width=\linewidth]{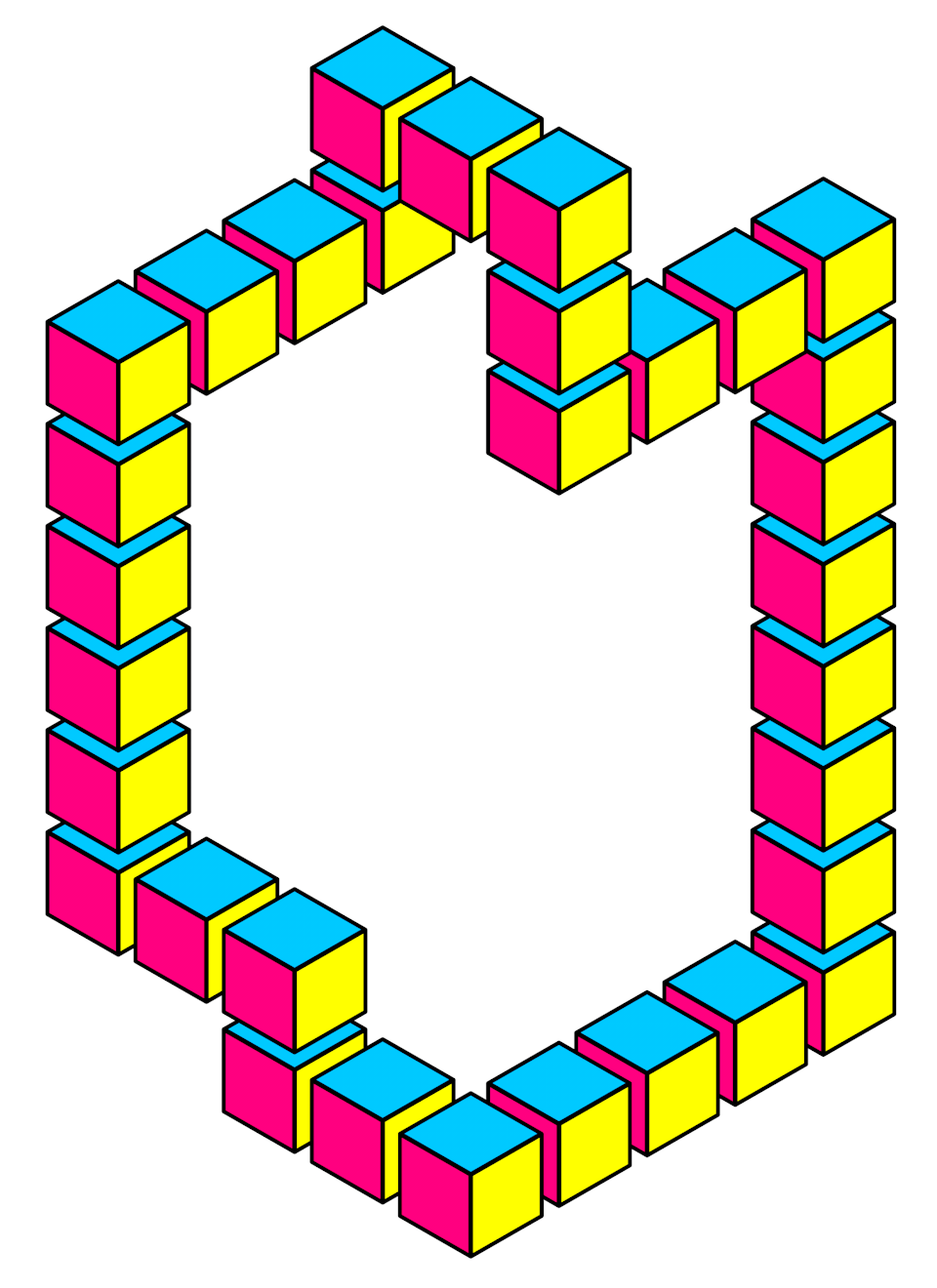}
    \end{minipage}
    \begin{minipage}[t]{0.24\textwidth}
        \centering
        \includegraphics[width=\linewidth]{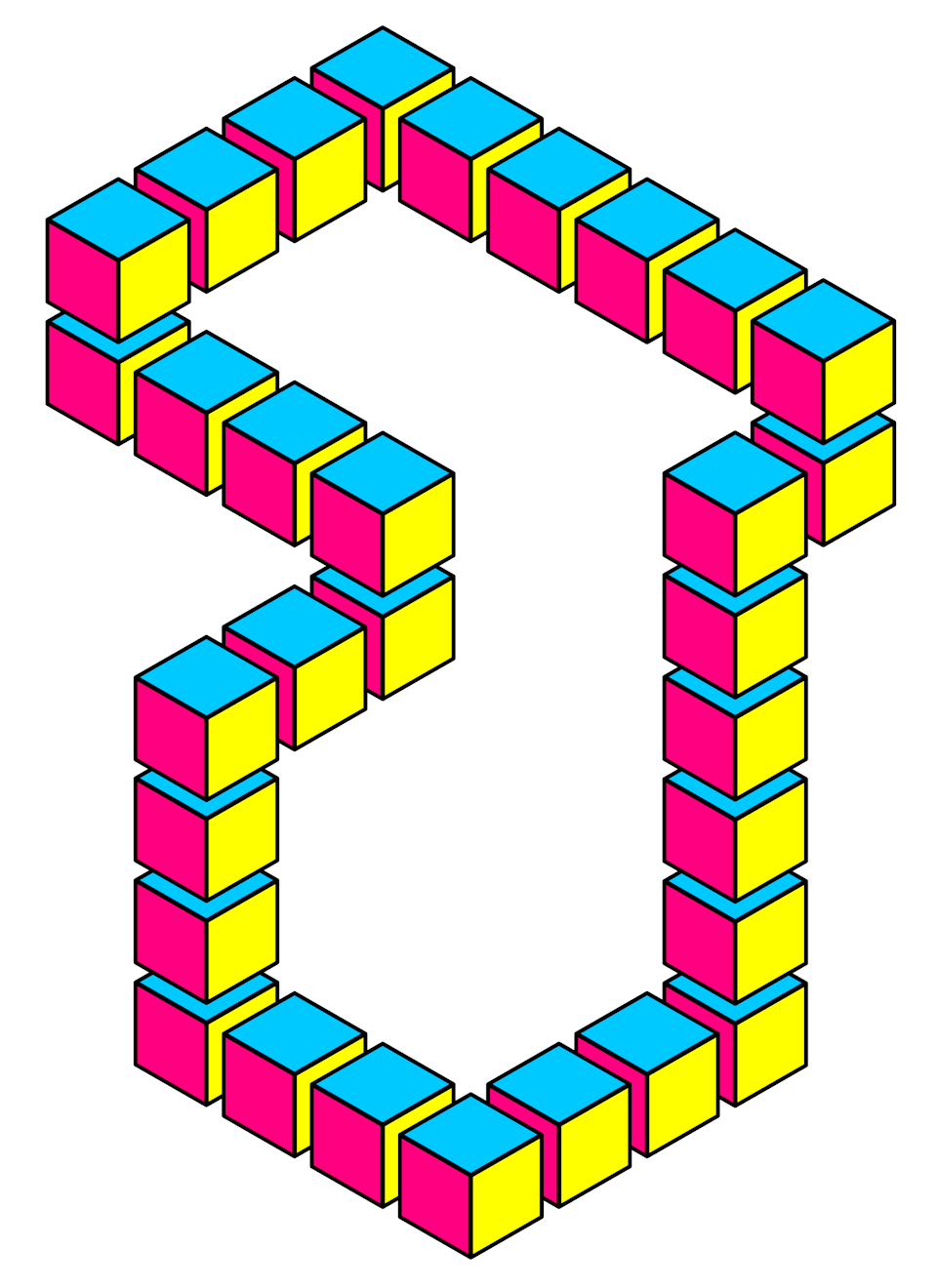}
    \end{minipage}
    \begin{minipage}[t]{0.24\textwidth}
        \centering
        \includegraphics[width=\linewidth]{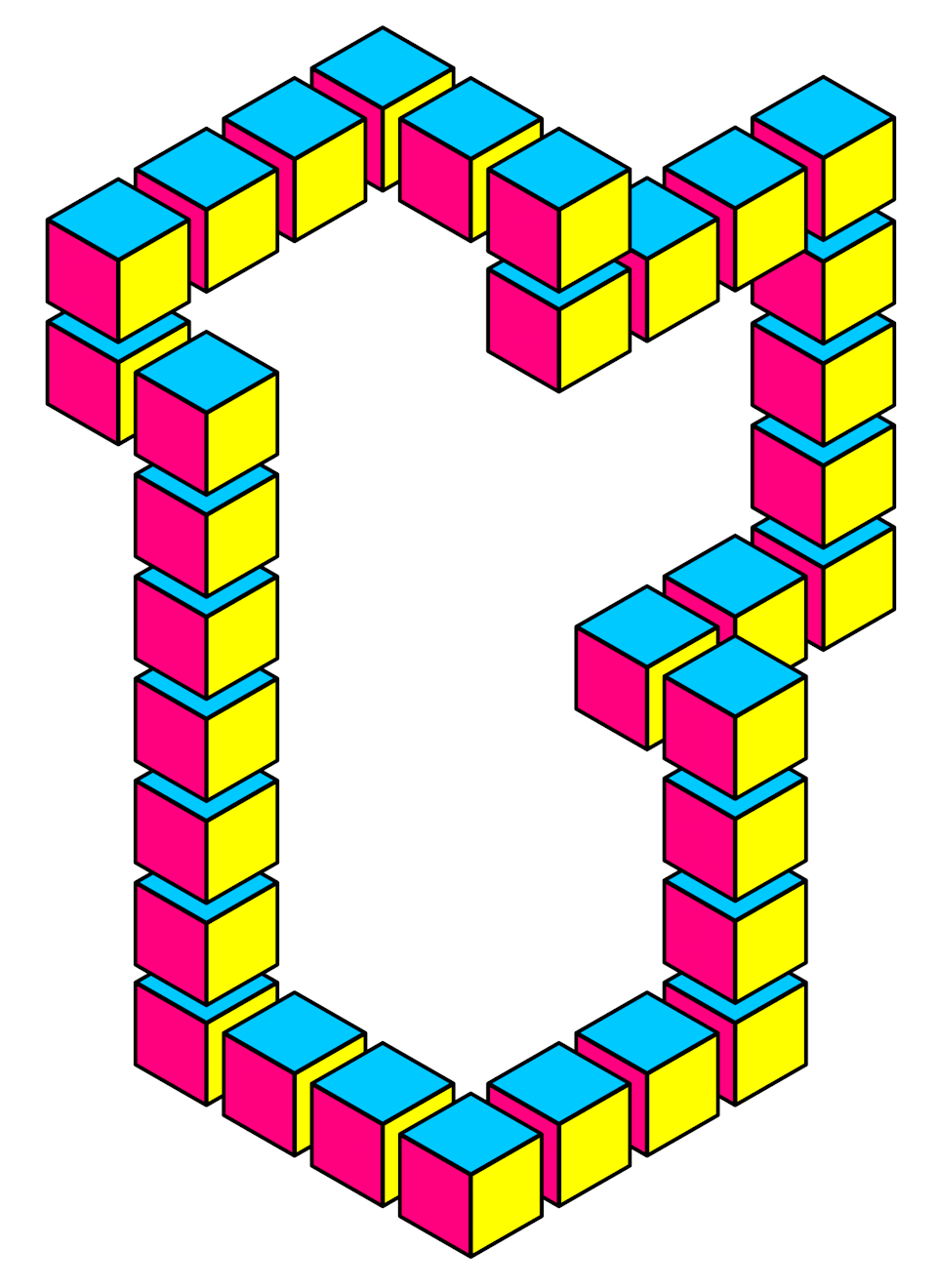}
    \end{minipage}
    \begin{minipage}[t]{0.24\textwidth}
        \centering
        \includegraphics[width=\linewidth]{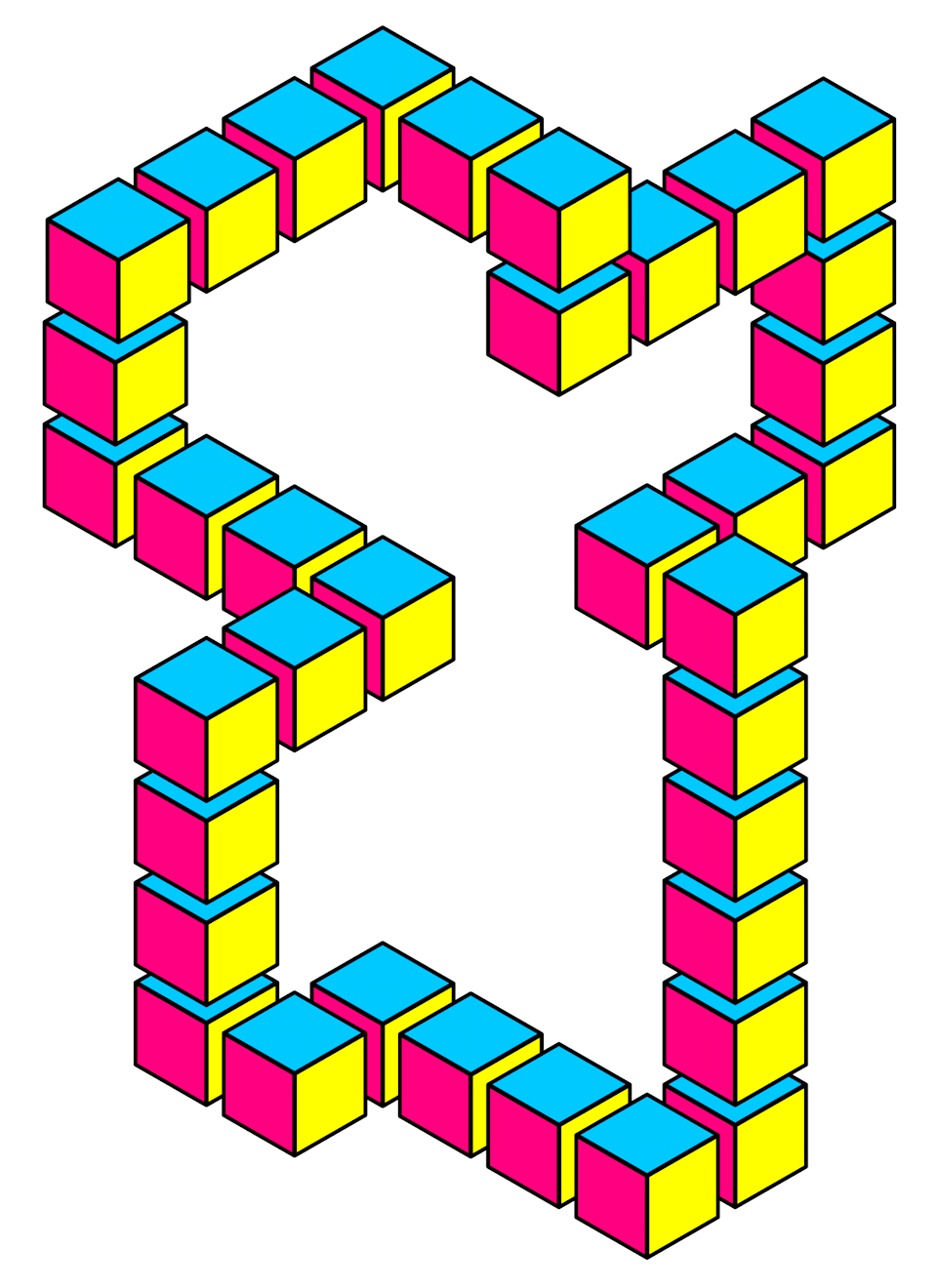}
    \end{minipage}
    
    \vspace{0.5em}
    
    \begin{minipage}[t]{0.24\textwidth}
        \centering
        \includegraphics[width=\linewidth]{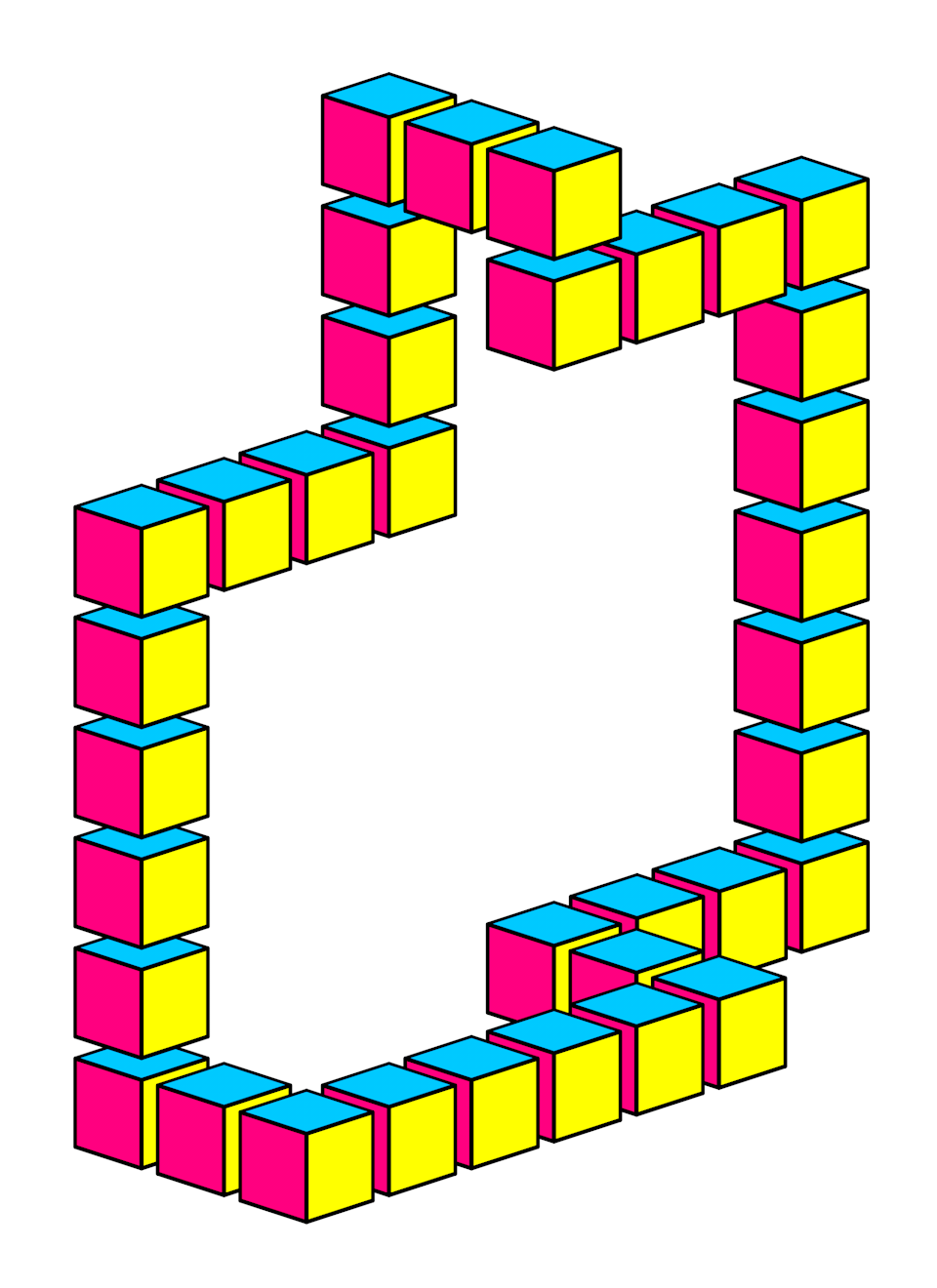}
    \end{minipage}
    \begin{minipage}[t]{0.24\textwidth}
        \centering
        \includegraphics[width=\linewidth]{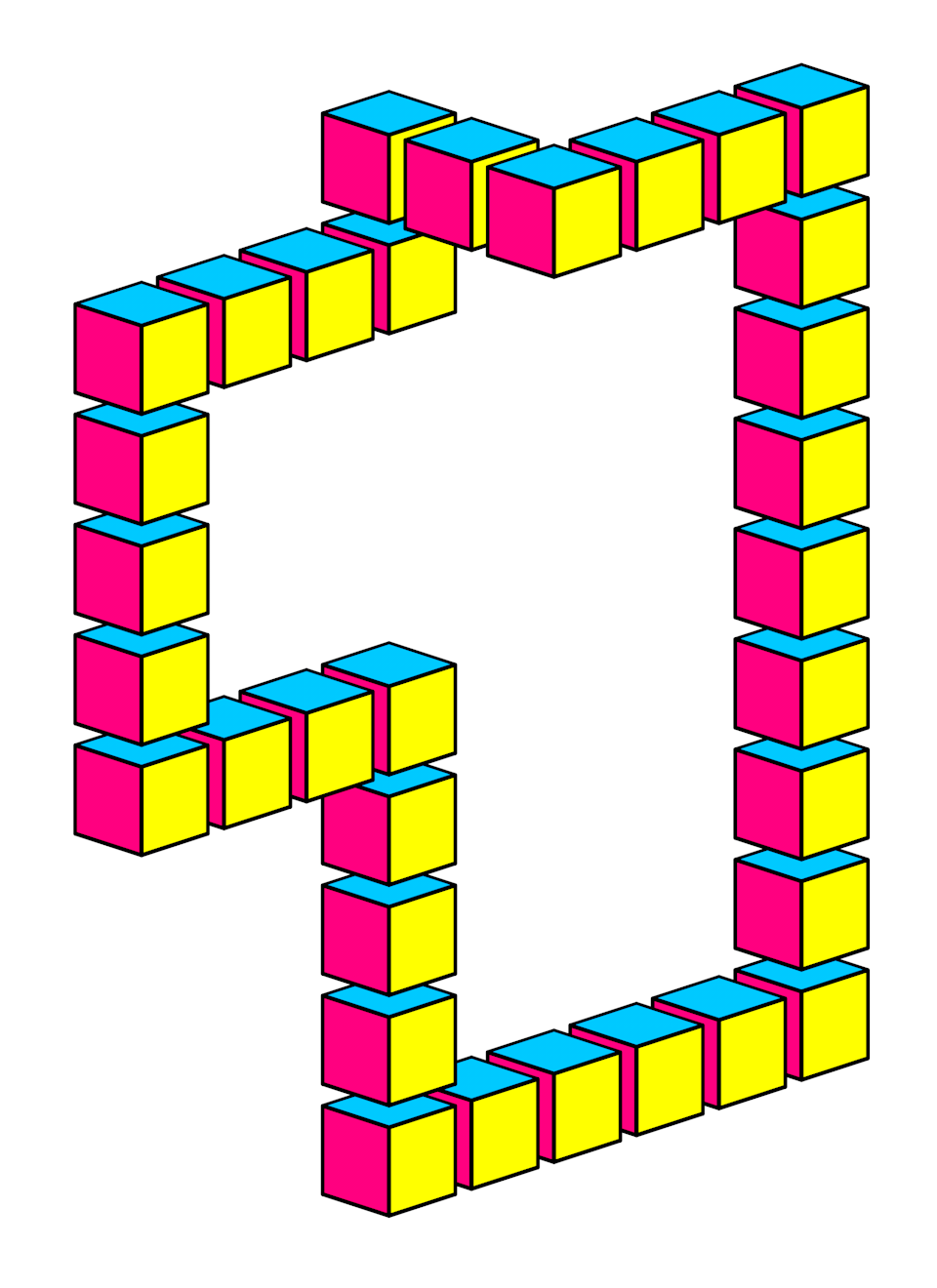}
    \end{minipage}
    \begin{minipage}[t]{0.24\textwidth}
        \centering
        \includegraphics[width=\linewidth]{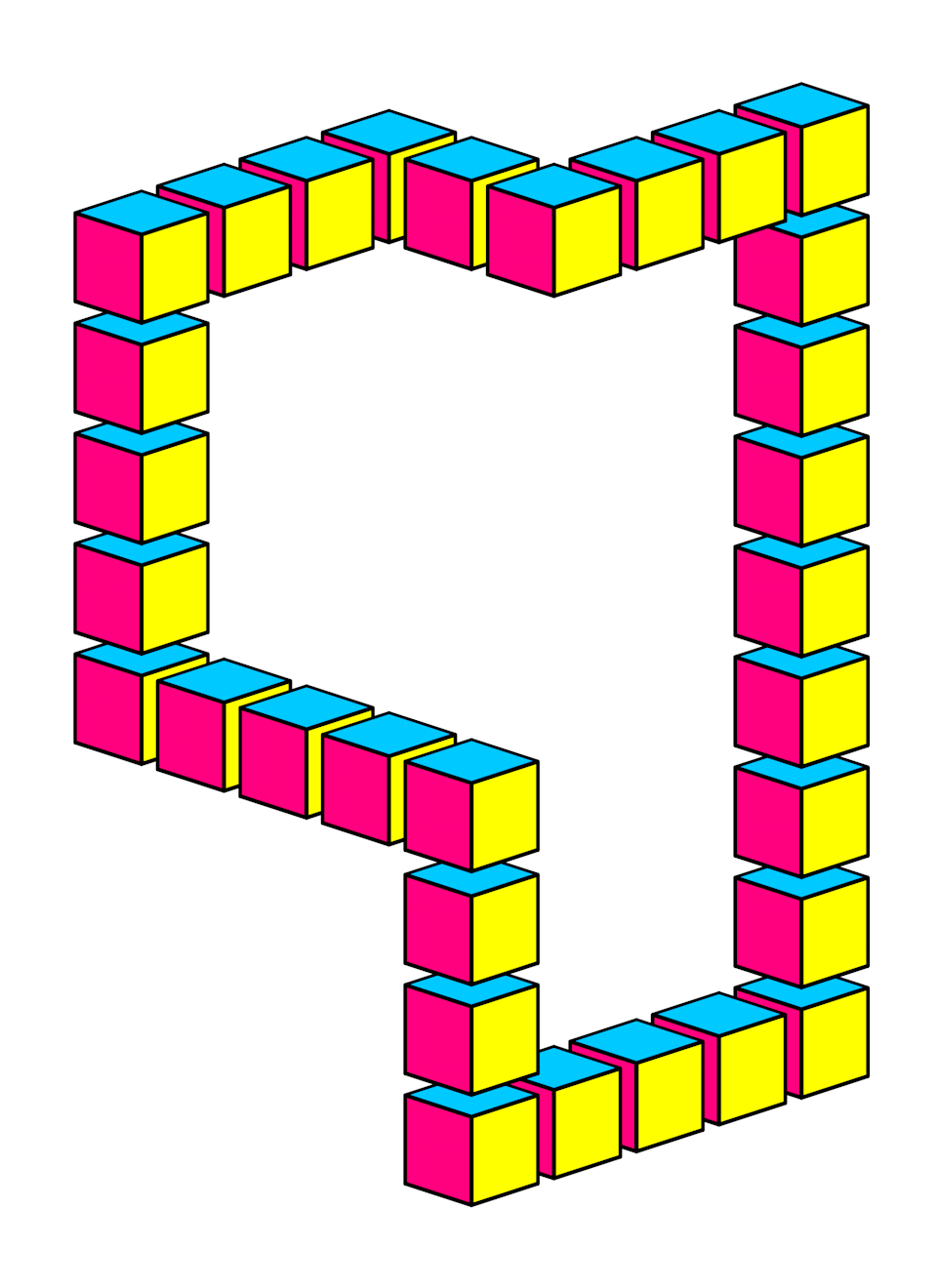}
    \end{minipage}
    \begin{minipage}[t]{0.24\textwidth}
        \centering
        \includegraphics[width=\linewidth]{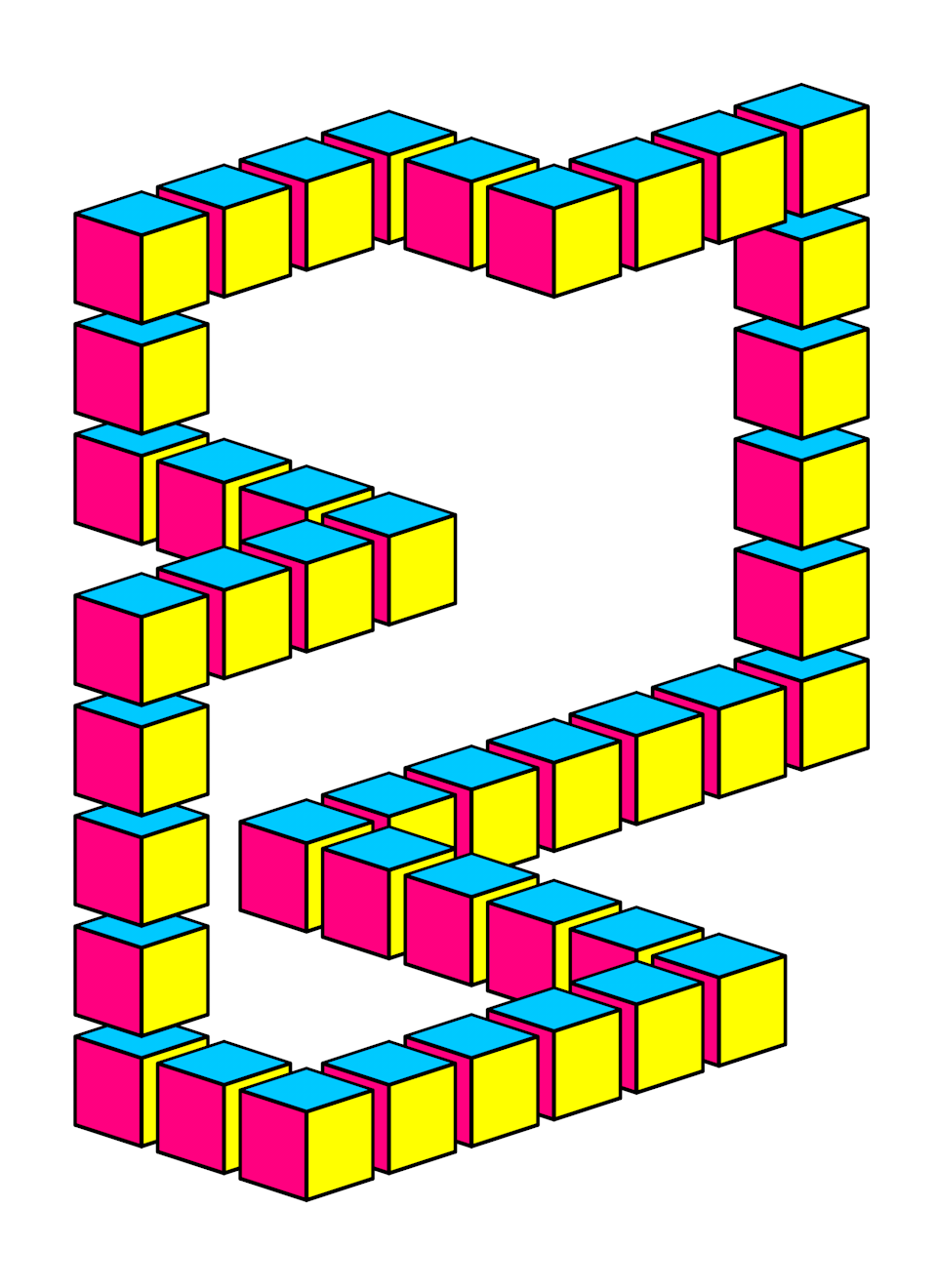}
    \end{minipage}
        
    \caption{\small Eight paradoxical loops in the style of Reutersv\"ard that exhibit local relative coordinate translations but do not possess corresponding global coordinates. These appear similar, but the precise mismatches of coordinates can vary. The top row is drawn in isometric projection, the bottom row in a dimetric projection.}
    \label{fig:penrose-stairs}
\end{figure}

Variants of the classic Penrose triangle and stairs are not uncommon in popular artworks. For example, Figure \ref{fig:penrose-stairs} contains eight examples of objects built from cubes in the style of Reutersv\"ard, with implicit unit translations between incident cubes in a local coordinate frame.  As one completes the circuit, there is a mismatch in the coordinates: a {\em holonomy}.

\begin{figure}
    \centering
    \setlength{\fboxsep}{0pt} % optional: no padding between box and image
    \fbox{\includegraphics[height=3.4cm]{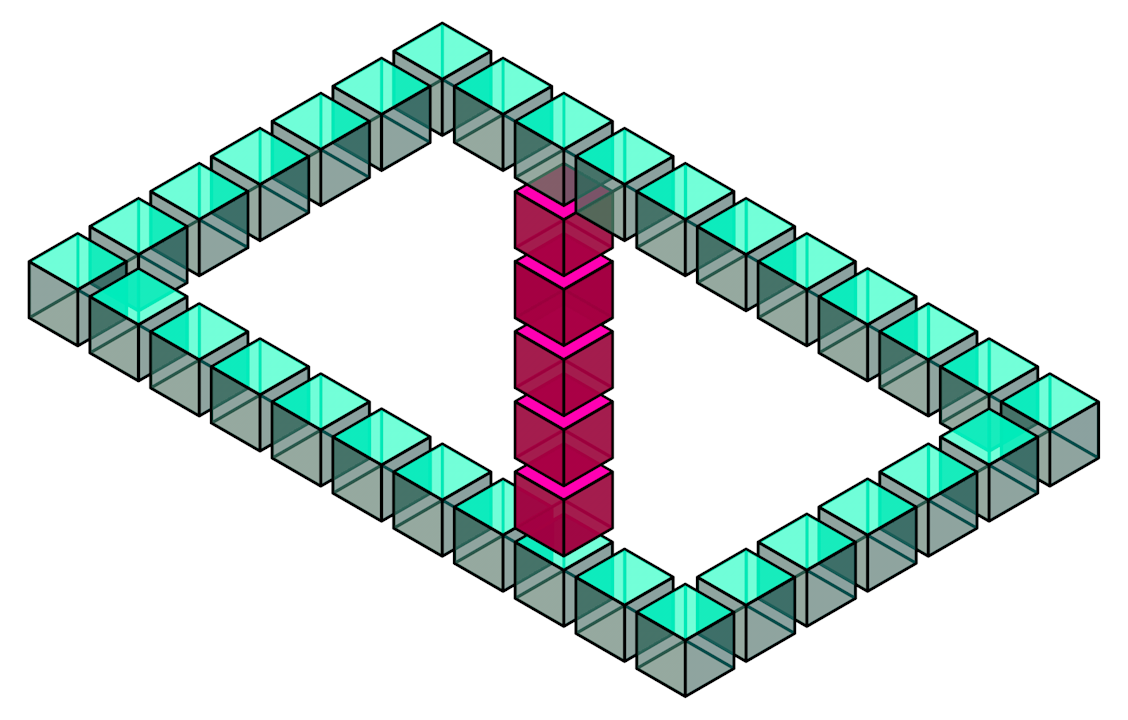}}%
    \hspace{0.5em}%
    \fbox{\includegraphics[height=3.4cm]{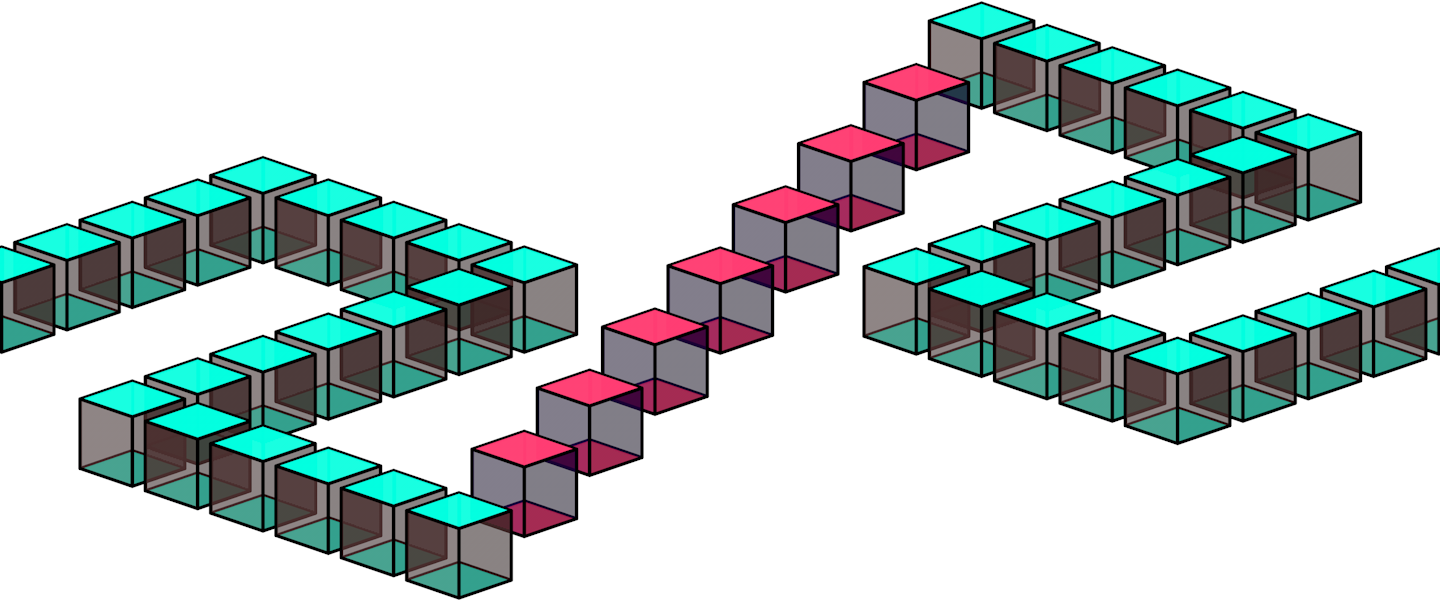}}
    \caption{\small [left] A variant of the classic Penrose triangle/staircase that has a ``ladder'' along which height changes. [right] A cylindrical version of an impossible stair: identify left and right sides to obtain a staircase which has well-defined height changes locally but not globally. Magenta steps connote height change; cyan steps are flat.}
    \label{fig:penrose_cyl}
\end{figure}

As another example, Figure \ref{fig:penrose_cyl}[left] is a variant that is meant to clearly show a paradox in height change.\footnote{Throughout, magenta steps will connote height change; cyan steps connote no height change.} One can travel the outer circuit without height change, but a shortcut up the vertical ``ladder'' leads to a paradoxical change in height. 

%------------------------------------------
\subsection{Cylindrical Staircase}
\label{sec:gallery-cyl}
%------------------------------------------

We move from loop-based paradoxes to those on identification spaces. Imagine a staircase that wraps around a cylinder. After completing one full circuit, we return to our starting point, yet the visual impression suggests we should be at a different height. The essential feature of a closed loop along which a quantity (say, height) undergoes a net non-zero change is present.

Such a cylinder can be realized as an identification space with two opposite sides identified: see Figure \ref{fig:penrose_cyl}[right]. Local geometric perceptions are preserved across this identification. As one crosses from the right edge to reappear at the left, the perceived height is to remain unchanged. This means that traversing the staircase involves a consistent increase in height when climbing the magenta stairs, with no height change along the cyan paths or when crossing the identification boundary. Following this path continuously around the cylinder results in returning to the starting point at a different height. 

\begin{figure}
    \centering
    \setlength{\fboxsep}{0pt} 
    \fbox{\includegraphics[width=4.5in]{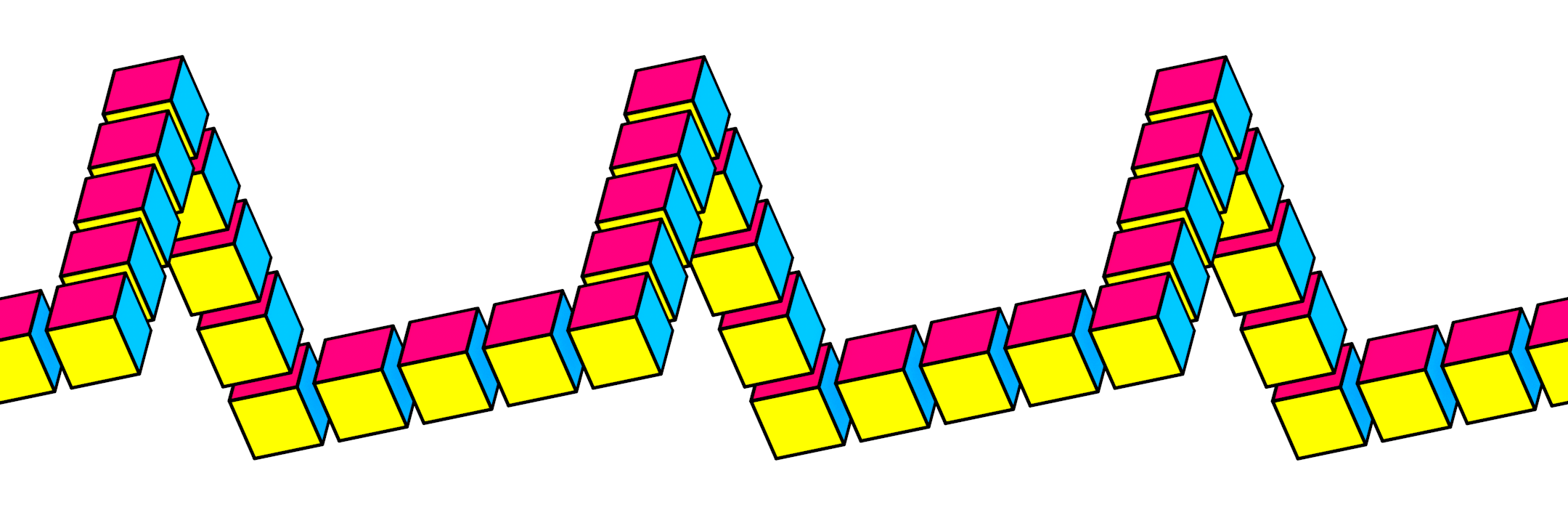}}
    \caption{\small An impossible object on a cylinder with left and right sides identified. How, precisely, should one characterize the holonomy?}
    \label{fig:long_cyl}
\end{figure}

A different presentation of a cylindrical Penrose-type stair appears in Figure \ref{fig:long_cyl}, which has left and right sides identified. The mismatch here is less about height and more similar to that of Figure \ref{fig:penrose-stairs}. One senses that this is a type of geometric ``covering'' of the classic Penrose triangle.

%------------------------------------------
\subsection{Orientation Flips: M\"obius Staircase}
\label{sec:gallery-mobius}
%------------------------------------------

\begin{figure}[htb]
    \centering
    \setlength{\fboxsep}{0pt} % optional: no padding between box and image
    \fbox{\includegraphics[height=5cm]{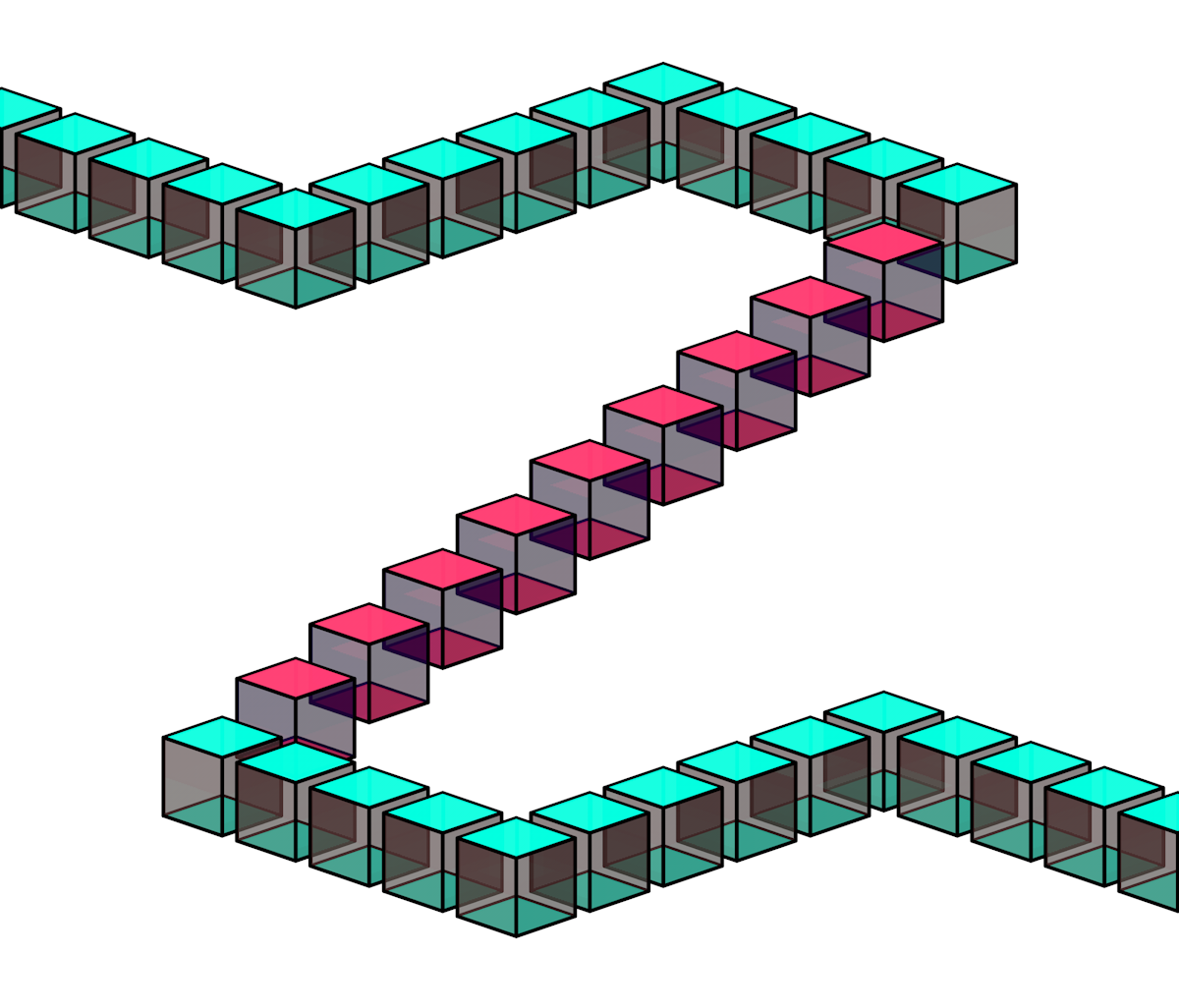}}%
    \hspace{3.5em}%
    \fbox{\includegraphics[height=5cm]{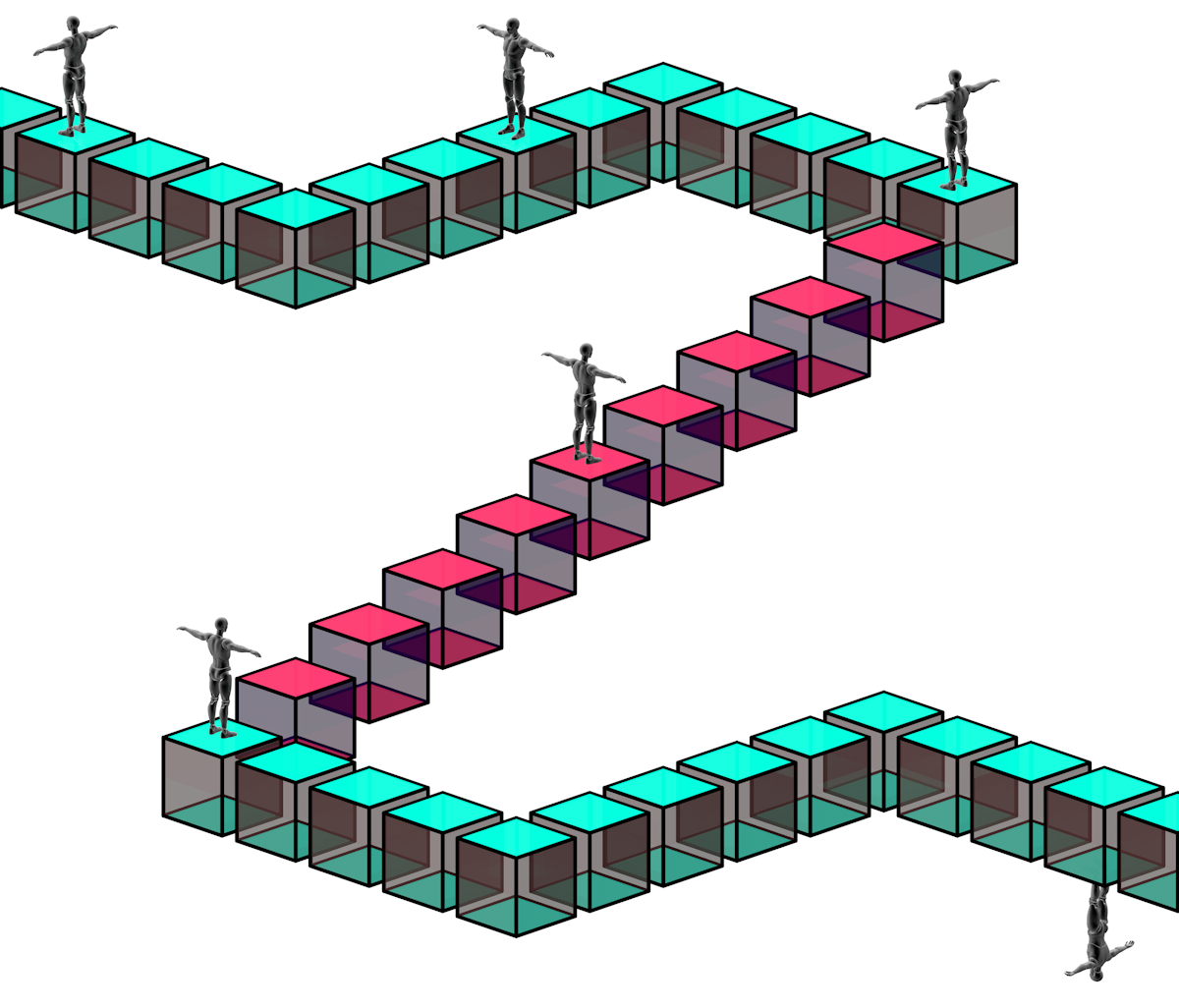}}
    \caption{\small [left] A paradoxical staircase on a M\"obius strip, with left and right identified with a vertical reflection. [right] This reflection reverses which ``side'' of the steps one walks on as one passes the left/right edge.}
    \label{fig:mobius}
\end{figure}

We can build a paradoxical staircase on a M\"obius strip as an identification space: see Figure \ref{fig:mobius}[left]. This, like the cylindrical case, represents the ambient space as a square with left and right sides identified; but now, the side edges are glued with a reflection. This reflection alters both the vertical orientation within the square and the perceived sense of height extending from the surface.

As hinted in Figure \ref{fig:mobius}[right], for the staircase as shown, someone beginning at the bottom of the steps walks up the stairs and continues to the left. Past the top left edge, they appear at the bottom right edge, but now walking on the ``underside'' of the surface. If they continue forward, then what they had experienced as ascending stairs the first time are now experienced as descending stairs on the second loop due to an inverted orientation of height.

\begin{remark}
\label{rem:orientation-preserving}
The reflection across the seam is a rule imposed by the model, not a consequence of physics. Crossing the seam, the local sense of \emph{up} reverses: the top of a tread is now seen from beneath. We record this as a reflection of the signed height measured along the local normal to the band. That coordinate is not ambient Euclidean height. The stair cues assert, at every tread, that the local up is the local normal. The physical argument needs one more assumption, which the drawing invites: that every local up is also the single ambient vertical. Under that assumption a physical band cannot exist: a normal line field that is everywhere vertical is trivial, and a band with trivial normal line field in $\mathbb{R}^3$ is orientable. The contradiction the figure presents lies exactly here, between the local up demanded by the stair cues, which reverses across the seam, and a global up, which cannot. For the projective plane, which does not embed in $\mathbb{R}^3$ at all, only the identification rule survives.
\end{remark}

\begin{figure}[hbt]
    \centering
    \setlength{\fboxsep}{0pt} % optional: no padding between box and image
    \fbox{\includegraphics[height=5.5cm]{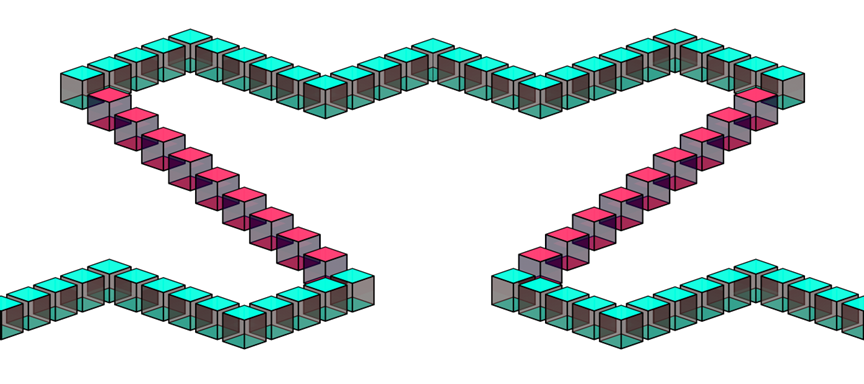}}
    \caption{\small The double cover of the M\"obius strip is an annulus, and the staircase of Figure \ref{fig:mobius} lifts to a staircase on a cylinder. A climber perceives walking along the M\"obius  staircase twice in the same direction and yet returns to their initial position, height, and orientation. Along the way, left-right and up-down have been switched when crossing the boundary of the fundamental domain.}
    \label{fig:mobiuscover}
\end{figure}

It is helpful to illustrate the double cover of the strip as an orientable annulus, as in Figure \ref{fig:mobiuscover}. From the perspective of a stair-climber, continual progression in their local frame leads to alternating left-right turns as perceived in the universal covering space. Traversing the loop once inverts orientation, while a second traversal restores the original orientation. After two complete circuits, the climber returns to the starting point with the same local orientation and height.

The holonomy is a reflection of order two. Starting at the foot of the flight, one circuit produces a height change of $\Delta$ together with an orientation flip, and a second circuit undoes both. A walker who starts at a different step sees a different one-loop change -- $\Delta - 2j$ from step $j$, so $\Delta$ from the foot, $-\Delta$ from the top, and $0$ from the middle step when $\Delta$ is even -- but always the same cancellation after two circuits. The paradox is torsion in a structure group that also contains translations of infinite order, and it does not rely on a bistable percept such as the Necker cube.

%------------------------------------------
\subsection{Climbing the Projective Plane}
\label{sec:gallery-RP2}
%------------------------------------------

\begin{figure}[hbt]
    \centering
    \setlength{\fboxsep}{0pt} % optional: no padding between box and image
    \fbox{\includegraphics[height=6cm]{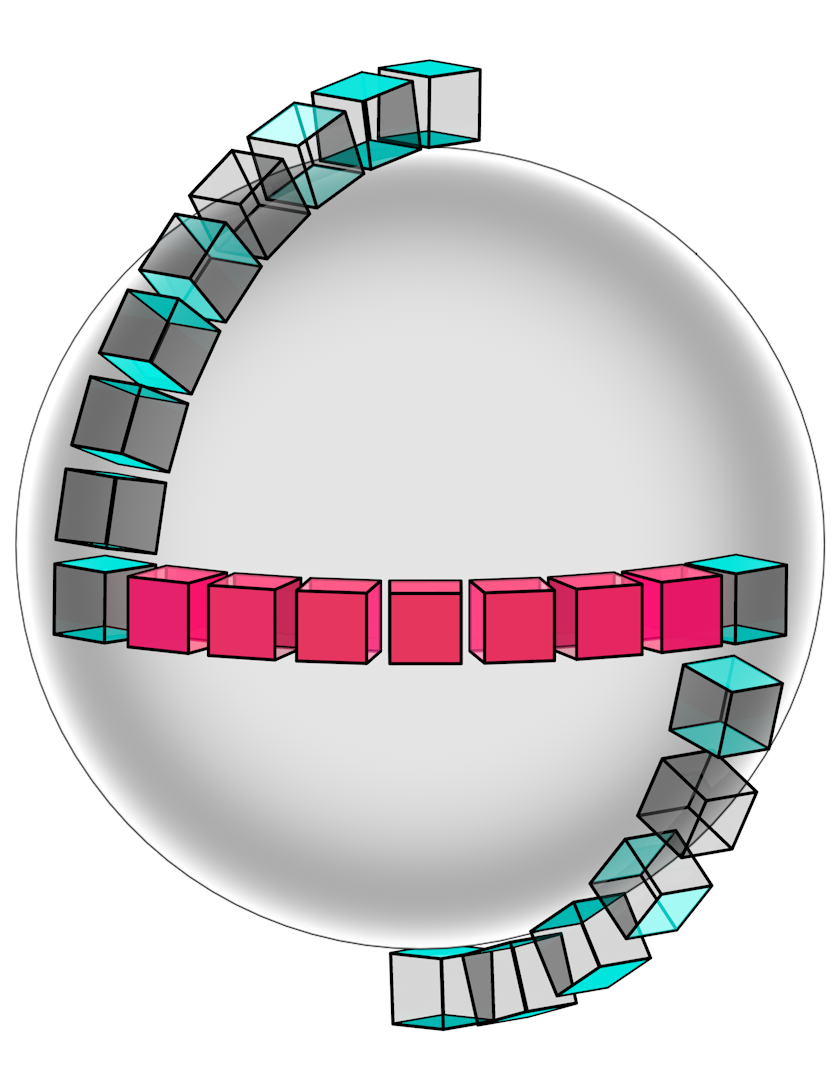}}%
    \hspace{0.5em}%
    \fbox{\includegraphics[height=6cm]{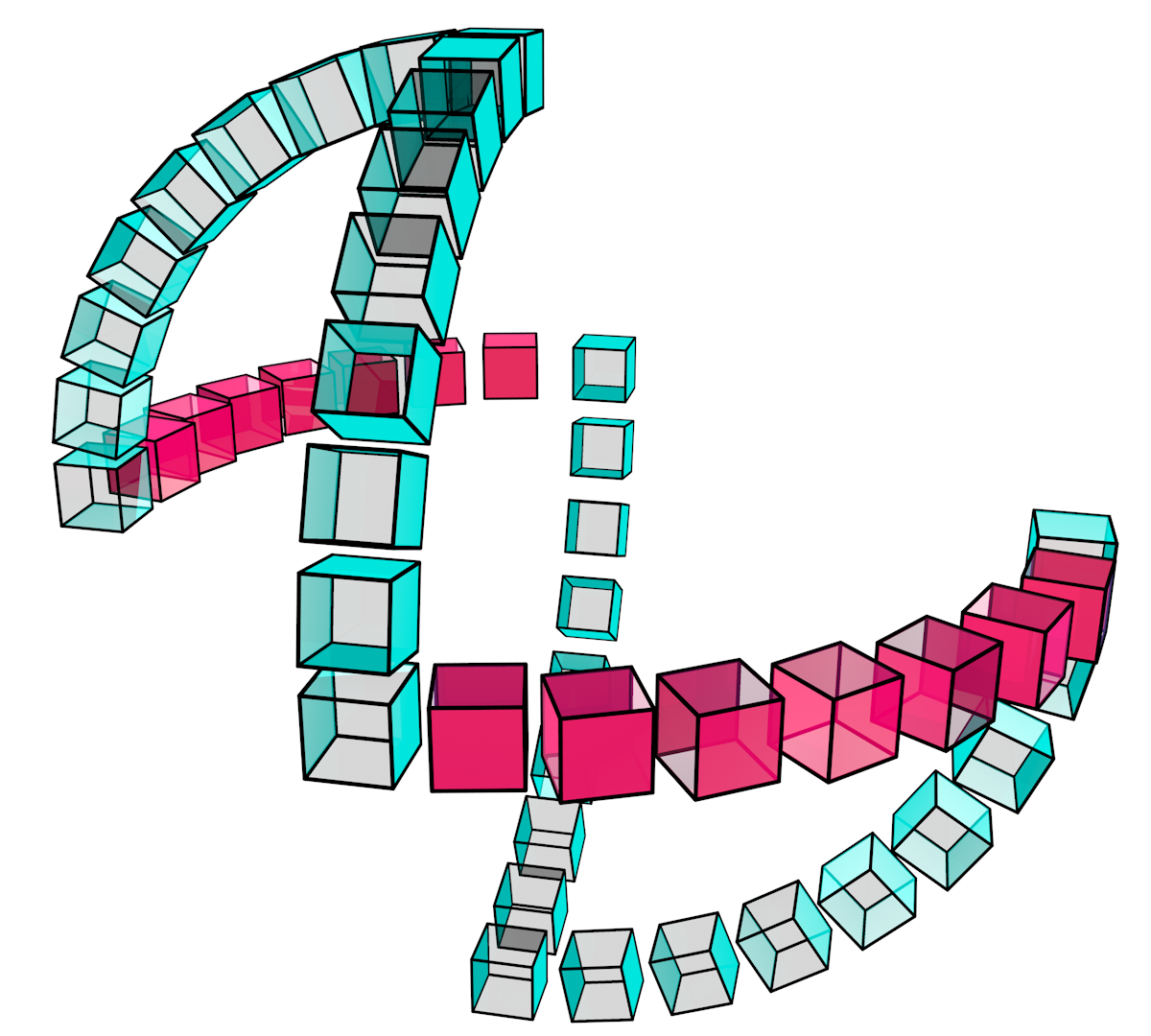}}%
\caption{\small Impossible staircase on a real projective plane $\mathbb{R}P^2$ as a hemisphere with antipodal boundary points identified with orientation reversal [left], where the north and south poles, which lie on the boundary circle, are identified with a reflection. Note the resemblance to the classical Penrose triangle -- three turns. The universal cover $S^2$ with the lifted staircase [right] demonstrates the torsional paradox. To illustrate the hemisphere [left], a parallel camera is used, which suppresses the stair-height change.}
\label{fig:RP^2}
\end{figure}

The real projective plane $\mathbb{R}P^2$ is a M\"obius band with a disc sewed to its boundary, and thus, the previous example extends naturally to a paradoxical staircase on the projective plane. However, a more satisfying example can be illustrated by fitting a staircase to a hemisphere with antipodal boundary points identified, then lifting to the double cover $S^2$: see Figure \ref{fig:RP^2}. 

Imagine someone beginning at the north pole of the sphere, walking to the equator with no height change. Then, after a left turn, the climber ascends the equatorial stairs. Turning to the right, the walker continues on a flat path to the south pole. The antipodal map inverts the climber back to the north pole, but with the radial coordinate also inverted, as prescribed by the identification rule of Remark \ref{rem:orientation-preserving}. The climber -- always alternating left/right turns -- walks another loop on $\mathbb{R}P^2$, but this time the stairs are perceived as descending (walking along the back side). Two loops, and the climber is back at the north pole at the original height.

This is the same phenomenon as with the M\"obius staircase: again, a $\mathbb{Z}_2$ torsion. It is interesting that the same loop in the M\"obius band and in $\mathbb{R}P^2$ have different behaviors over the ambient space's fundamental group; yet the paradox's nature is the same.

%%%%%%%%%%%%%%%%%%%%%%%%%%%%%%%%%%%%%%%%%%%%%%%%%%
%%%%%%%%%%%%%%%%%%%%%%%%%%%%%%%%%%%%%%%%%%%%%%%%%%

%------------------------------------------
\subsection{Zigzag Depth Paradox}
\label{sec:gallery-zigzag}
%------------------------------------------

\begin{figure}[hbt]
    \centering
    \setlength{\fboxsep}{0pt} \fbox{\includegraphics[width=5.5in]{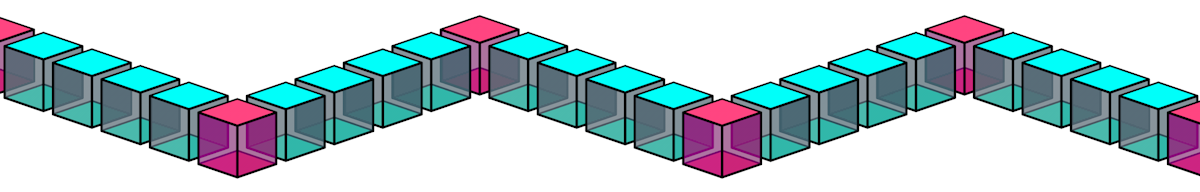}}%
    
    \fbox{\includegraphics[width=5.5in]{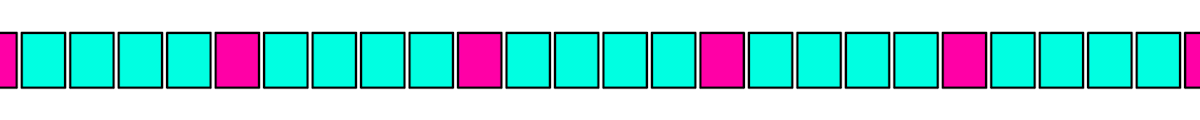}}%
\caption{\small Zigzag depth paradox: A path with alternating right-angled turns [top] can be realized as a loop (identifying left and right sides) only when the number of corner points is even. When viewed from the side marking only straight-vs-corner points [bottom], such a loop is sensible for any number of corner points, but an odd number presents a binary visual paradox that does not lift to consistent alternating depths.}
\label{fig:zigzag}
\end{figure}

A distinct class of visual paradoxes emerges when we consider depth ambiguity rather than height differentiation. Consider the zigzag path constructed of segments which alternates direction at each corner in Figure \ref{fig:zigzag}[top]. When extended infinitely on each side or wrapped around a cylinder there is no paradox.

%%%%%%%%%%%%%%%%%%%%%%%%%%%%%%%%%%%%%%%%
%%%%%%%%%%%%%%%%%%%%%%%%%%%%%%%%%%%%%%%%

Imagine now viewing this system from the side, illustrated in a manner that records only the corner points without indicating the depth, as in Figure \ref{fig:zigzag}[bottom]. This yields a relative measurement: the viewer knows that the path is changing direction but without knowing which way.

In isolation, such a path presents mere bistable perception, similar to the Necker cube effect. However, when its endpoints are identified to form a closed loop, a genuine paradox can emerge. With an odd number of corners, the resulting perceptual inconsistency is a binary ($\mathbb{Z}_2$-type) paradox based on depth perception. Unlike the M\"obius staircase which involves coupled height and orientation changes, the zigzag paradox focuses solely on the impossibility of consistently assigning one of two depth states (e.g., `protruding' vs. `receding') around a loop with an odd number of alternations. It is incidentally an example of a 1-dimensional projection that does not lift to a consistent 2-dimensional structure -- much like a Penrose triangle is a phenomenon of one dimension higher.

%------------------------------------------
\subsection{Commuting Paradoxes: Torus Staircase \& Ladder}
\label{sec:gallery-torus}
%------------------------------------------

\begin{figure}
    \centering
    \setlength{\fboxsep}{0pt}
$
\vcenter{\hbox{\fbox{\includegraphics[height=5cm]{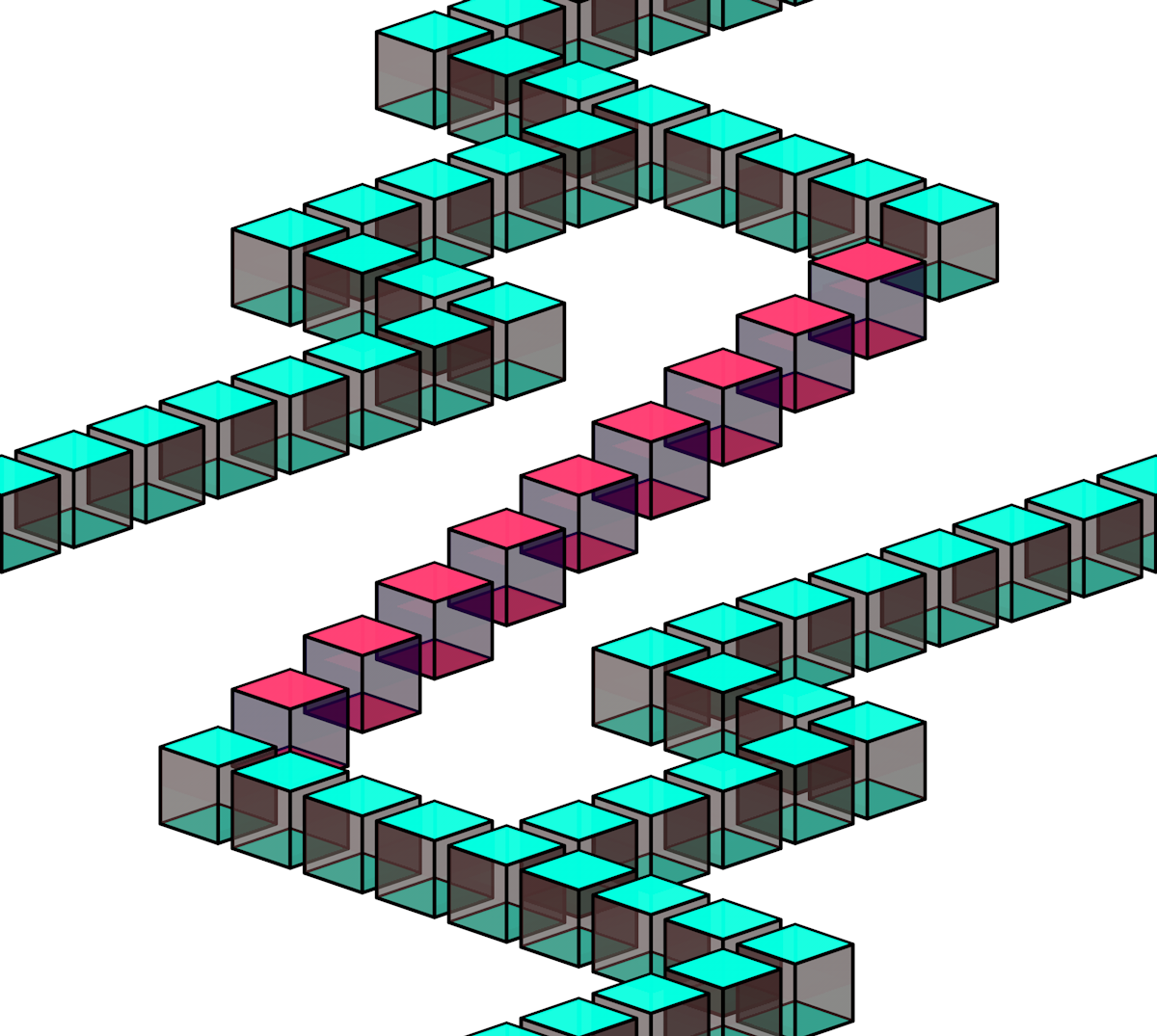}}}}
$
    \hspace{0.5em}
$
\vcenter{\hbox{\fbox{\includegraphics[height=8.5cm]{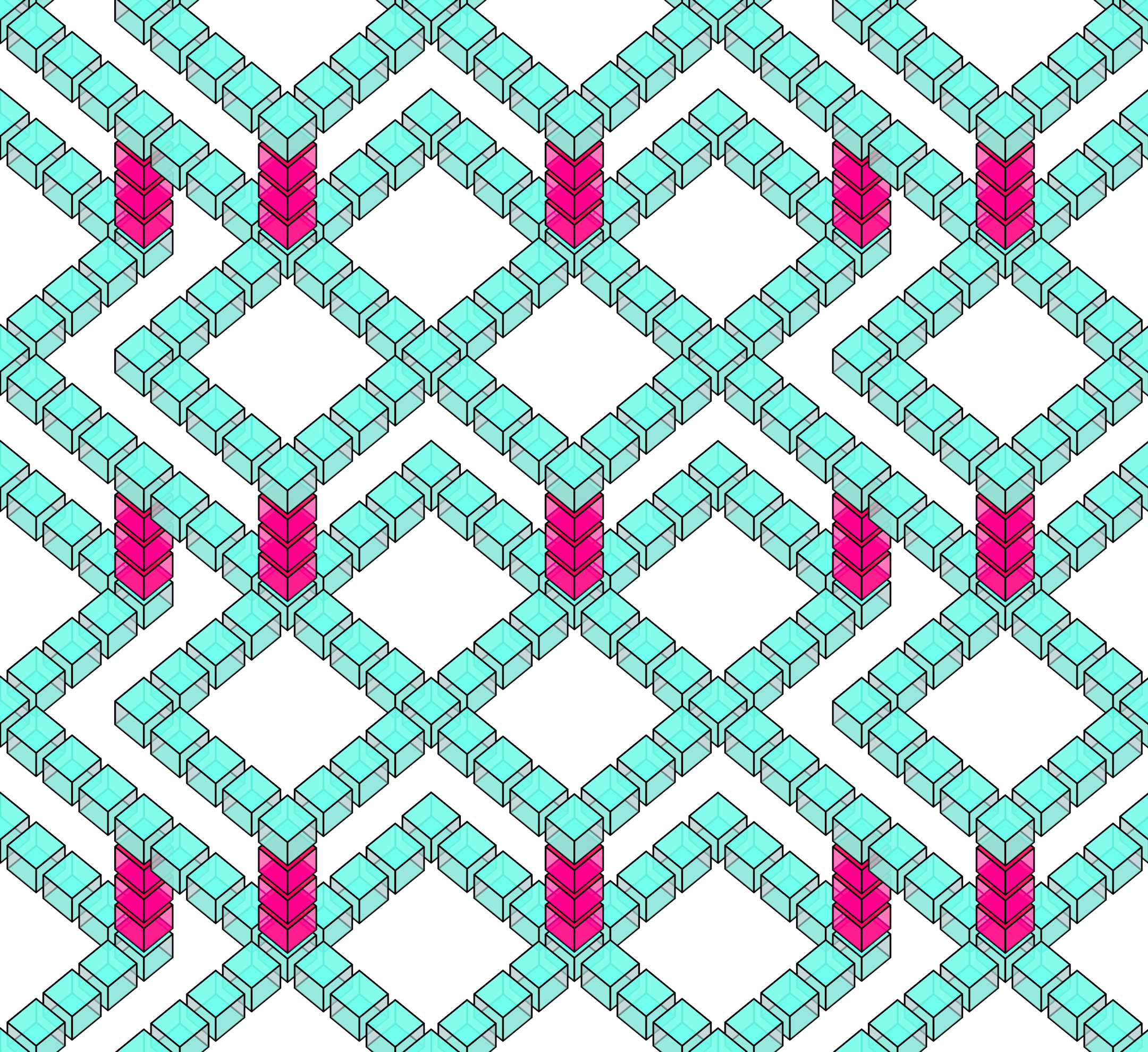}}}}
$
\caption{\small Impossible staircase on a torus $T^2$ represented as a square with opposite edges identified: [left] a simple example with a single height change; [right] a complex example from a different perspective, where a unit height change is indicated by vertical ladders in magenta. In this latter case, a simple horizontal loop induces two units of height change, whereas a simple vertical loop induces three units of height change.}
\label{fig:torus}
\end{figure}

A torus has two fundamental, independent loops (longitude and meridian) which commute in the fundamental group. Imagine a staircase pattern on a torus where traversing these loops results in a net height change. Figure \ref{fig:torus} gives examples on the identification space formed by a square with opposite sides identified. One notes that paths which remain within the square are not paradoxical. In Figure \ref{fig:torus}[left] the horizontal and vertical generating loops of the torus yield a height change of $+1$ each (up to units/orientation). For Figure \ref{fig:torus}[right], a simple vertical loop induces three units of height change; a simple horizontal loop induces two units of height change. 

%%%%%%%%%%%%%%%%%%%%%%%%%%%%%%%%%%%%%%%%%%%%%%%%%%
%%%%%%%%%%%%%%%%%%%%%%%%%%%%%%%%%%%%%%%%%%%%%%%%%%

These paradoxes, like the fundamental group of the torus, are abelian, as one can check by traversing a path along a commutator and examining the resulting height change.  

%------------------------------------------
\subsection{Nonabelian Paradox on the Klein Bottle}
\label{sec:gallery-klein}
%------------------------------------------

Taking the identification space for the torus and adding a twist to one pair of edges yields the Klein bottle $K$. Let $a$ be the loop that crosses the twisted pair of edges and $b$ the loop that crosses the untwisted pair. Then
\[
\pi_1(K) = \langle a, b \mid aba^{-1}b = 1 \rangle \cong \mathbb{Z} \rtimes \mathbb{Z},
\]
a nonabelian group in which conjugation by $a$ inverts $b$ and both generators have infinite order.

\begin{figure}[h]
    \centering
    \setlength{\fboxsep}{0pt} % optional: no padding between box and image
\fbox{\includegraphics[width=0.6\textwidth]{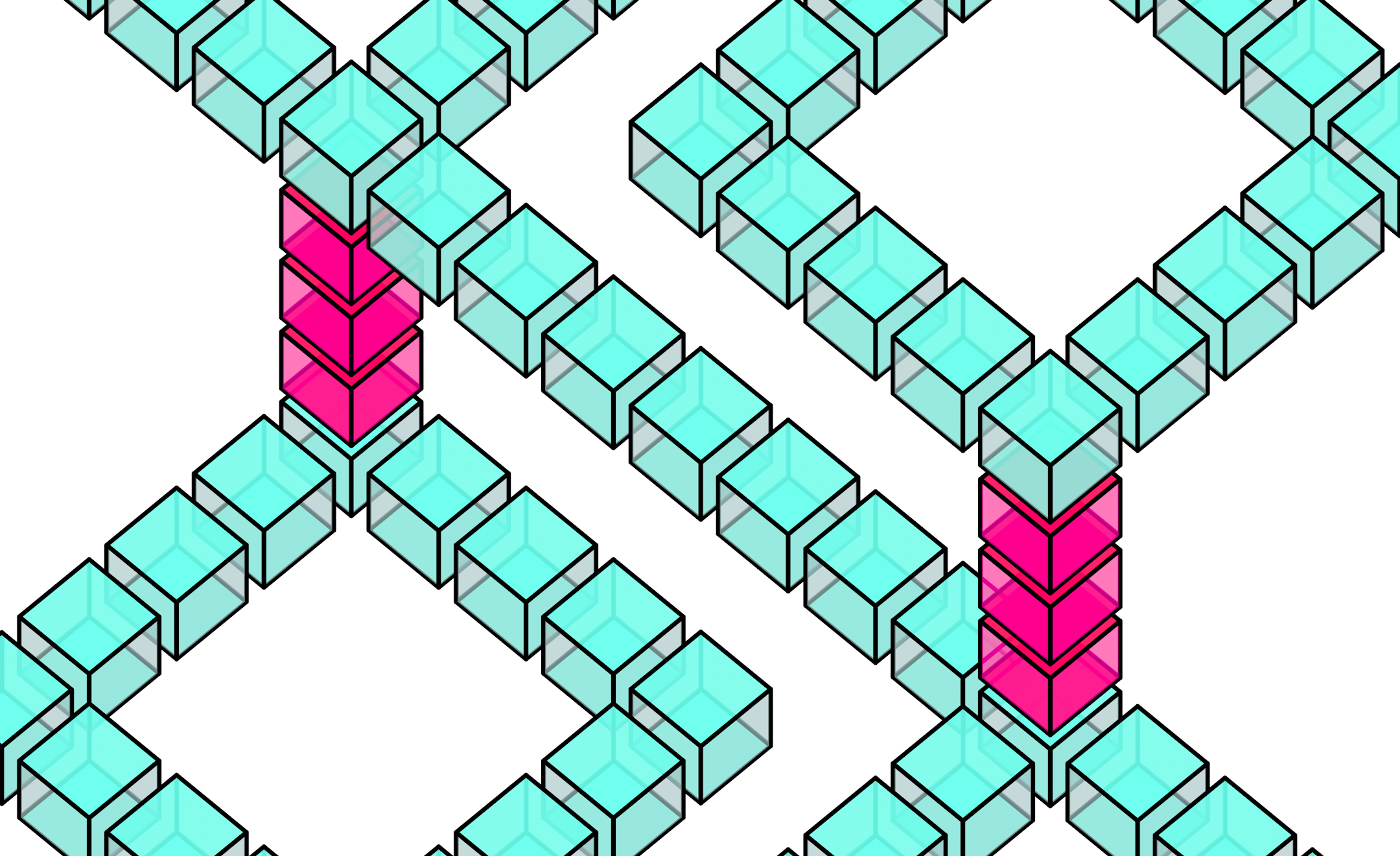}}

    \vspace{0.5em}
    \fbox{\includegraphics[width=\textwidth]{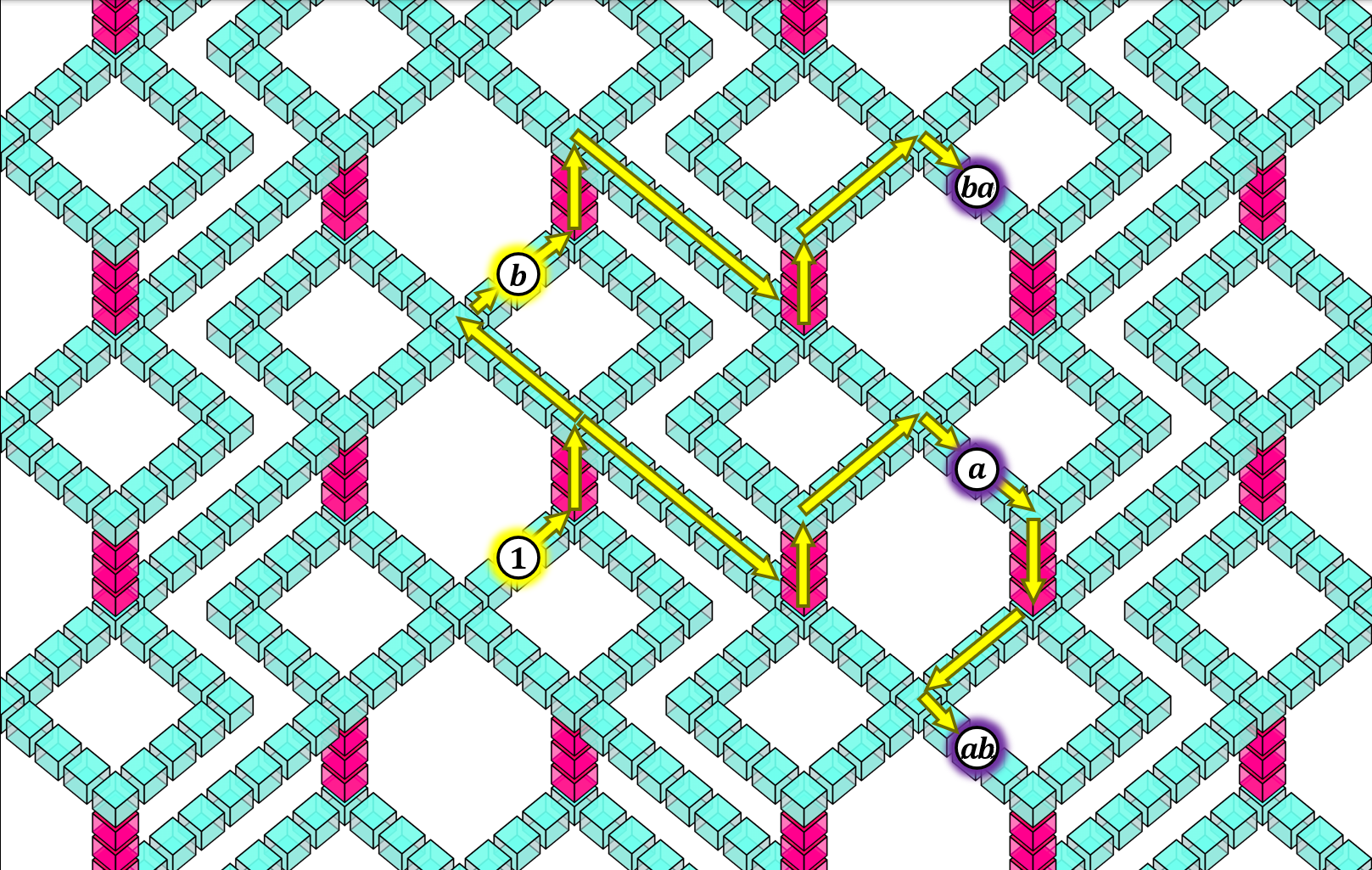}}
    \caption{\small The Klein ladder. [top] A fundamental domain for the Klein bottle, with top and bottom identified by translation and left and right identified with a vertical reflection; magenta ladders are vertical stacks of cubes, each climbing one unit of height. [bottom] A lift of the Klein ladder to the universal cover of the Klein bottle, drawn perceptually. From the basepoint $1$, the loop $b$ climbs one ladder; the loop $a$ climbs two ladders and arrives mirror-reflected, so that a subsequent $b$ descends. The words $ab$ ($a$ then $b$) and $ba$ ($b$ then $a$) end at different lifts of the basepoint, one and three ladders up: a noncommutative paradox.}
    \label{fig:klein}
\end{figure}

A \emph{Klein ladder} captures this relation. As in the torus ladder of Figure \ref{fig:torus}[right], chains of cubes run level, and height changes only on the ladders, vertical stacks of magenta cubes, each climbing one unit. Figure \ref{fig:klein}[top] is a fundamental domain: its left and right sides are identified with a vertical reflection and its top and bottom by translation, so that $a$ runs horizontally and $b$ vertically. Crossing the reflected seam turns every ladder upside down: a ladder climbed before the crossing is descended after it, the $\mathbb{Z}_2$ behavior of the M\"obius staircase. The universal cover, Figure \ref{fig:klein}[bottom], makes the interaction of the two loops visible. From the basepoint marked $1$, the loop $b$ climbs one ladder and the loop $a$ climbs two and arrives mirror-reflected, so that $b$ traversed after $a$ descends where $b$ traversed alone climbs. The words $ab$ and $ba$, read in the order of traversal, therefore end at different lifts of the basepoint, one ladder and three ladders up. The relator $aba^{-1}b$ closes up; the commutator $aba^{-1}b^{-1}$, which would close up on a torus, does not. The ladder does not see all of $\pi_1(K)$: its holonomy along $a$ is an involution, so the representation factors through the infinite dihedral group $D_\infty \cong \mathbb{Z} \rtimes \mathbb{Z}_2$, the quotient of $\pi_1(K)$ by $a^2 = 1$ (Section \ref{sec:nonabelian_paradox}). 

%%%%%%%%%%%%%%%%%%%%%%%%%%%%%%%%%%%%%%%%%%%%%%%%%%
%%%%%%%%%%%%%%%%%%%%%%%%%%%%%%%%%%%%%%%%%%%%%%%%%%

The holonomies along $a$ and $b$ do not commute: the figure is, to our knowledge, the first visual paradox with nonabelian holonomy.

%------------------------------------------
\subsection{Boundary-Induced Impossibility}
\label{sec:gallery-necker}
%------------------------------------------

\begin{figure}[h]
    \centering
    \setlength{\fboxsep}{0pt} % optional: no padding between box and image
\fbox{\includegraphics[width=6in]{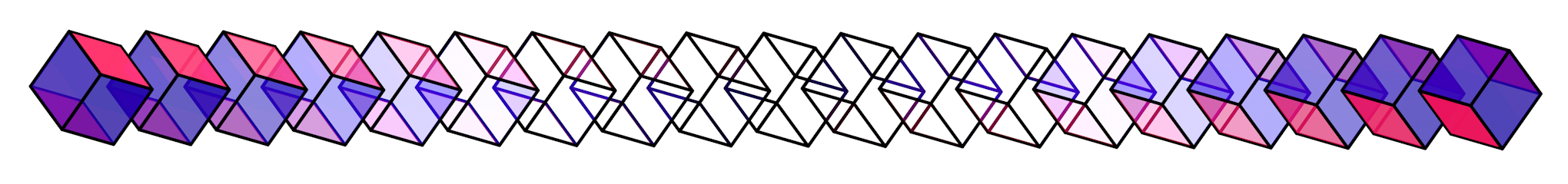}}%
  
\caption{\small The gradient Necker paradox is a sequence of cubes with forced contradictory interpretations at the endpoints, so that no consistent interpretation of the whole sequence exists.}
\label{fig:necker_gradient}
\end{figure}

Not all paradoxes arise from loops. Consider a sequence of {\em Necker cubes} (from Figure \ref{fig:classics3}[left]). For each cube, a viewer resolves the ambiguity by selecting one of two perceptions, propagating this choice locally. One can construct a paradoxical configuration where the cubes at opposite ends of the sequence are rendered with visual cues (shading and color) that force contradictory interpretations, as shown in Figure \ref{fig:necker_gradient}. The cubes gradually transition to wireframe representations, removing the visual cues that force a particular interpretation. This creates a region of ambiguity in the middle of the sequence and visually contradictory boundary constraints.

\begin{figure}[h]
    \centering
    \setlength{\fboxsep}{0pt} % optional: no padding between box and image
    \fbox{\includegraphics[width=6in]{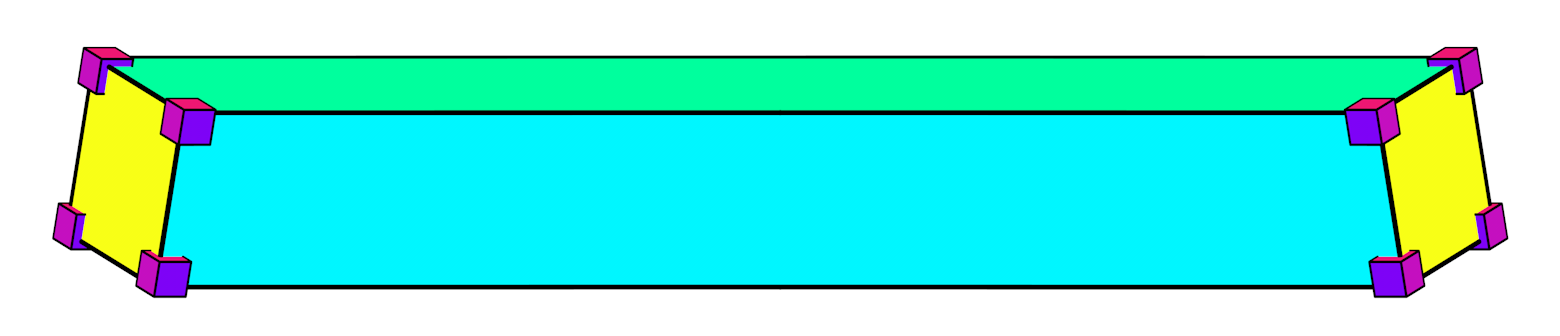}}%

    \fbox{\includegraphics[width=6in]{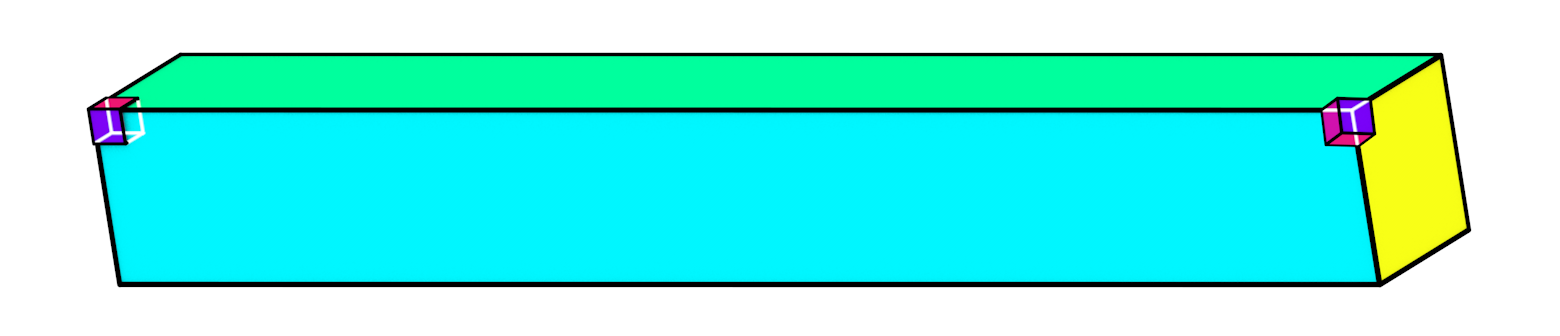}}%  
\caption{\small Impossible Bars: Two long bars with forced contradictory perspectives indicated via right-angle cues at the ends. The interior of the bar suggests continuity, causing a conflict with the endpoints. Covering one end of the bar then the other triggers a gestalt switch. The first bar [top] uses a contradictory perspective framing of the ends; the second [bottom] reverses interior and exterior corners. }
\label{fig:boundary}
\end{figure}

There are many such ``impossible objects'' that appear in popular artworks. One has a pair of conflicting forced interpretations at the ends of an interval, with a notion of continuation across the interval that makes the illusion stick. Figure \ref{fig:boundary} gives two long {\em Impossible Boxes} with perspectives forced by visual cues. 

% -----------------------------------------
\subsection{Observations and Questions}
% -----------------------------------------

This gallery of classical and novel visual paradoxes reveals several distinct types of impossibility. Some arise from simple loops (Penrose), others from non-orientable surfaces (M\"obius, $\mathbb{R}P^2$), commuting loops (torus), non-commuting interactions (Klein ladder), or boundary conditions (Necker steps, impossible bars). The resulting visual effects include monotonic ascent/descent, orientation flips, bistable states, and more.

The sections that follow give each figure an explicit model: a base, a structure group, a set of local states, and transport. The toolkit of sheaves, cohomology, and torsors comes first.

%%%%%%%%%%%%%%%%%%%%%%%%%%%%%%%%%%%%%%%%%%%%%
\section{Sheaves \& Sheaf Cohomology}
\label{sec:toolkit}
%%%%%%%%%%%%%%%%%%%%%%%%%%%%%%%%%%%%%%%%%%%%%

The visual paradoxes introduced in Section \ref{sec:gallery} arise from a fundamental tension between local consistency and global impossibility. To analyze this phenomenon rigorously, we employ three interlocking mathematical tools: sheaf theory, cohomology, and torsors. Sheaves provide the framework for organizing locally defined data on a space, capturing how geometric information (such as heights or orientations) exists at each point of a figure. Sheaf cohomology measures the obstruction to globalizing this local data. Torsors are the geometric objects that manifest this obstruction -- spaces where relative relationships are well-defined but absolute values cannot be consistently assigned. 

All three of these tools exist in the literature on sheaf theory \cite{giraud1971cohomologie,kashiwara1990sheaves,maclane1992sheaves}. Rather than wield the full imposing structure, we discretize to the combinatorial setting of a network base space -- a graph of vertices and edges. Network sheaves and their cohomology have recently proven useful in consensus \cite{riess2022diffusion}, opinion dynamics \cite{hansen2021opinion}, graphic statics \cite{cooperband2023homological,cooperband2023cosheaf}, neural networks \cite{hansen2020sheaf}, and more. 

%------------------------------------------
\subsection{Network Sheaves: Modeling Local Data on Graphs}
%------------------------------------------

Many impossible figures possess an underlying structure well-approximated by a graph $X = (V, E)$, consisting of vertices $V$ and edges $E$. Network sheaves provide a natural way to encode locally defined data on such 1-dimensional structures.

\begin{definition}
Let $X = (V, E)$ be a graph with an orientation on each edge. A \emph{network sheaf} $\mathcal{F}$ on $X$ with values in a category $\mathbf{C}$ assigns:
\begin{itemize}
\item To each vertex $v \in V$, an object $\mathcal{F}(v) \in \mathbf{C}$, called the \emph{vertex stalk} of $\mathcal{F}$ at $v$.
\item To each edge $e \in E$, an object $\mathcal{F}(e) \in \mathbf{C}$, called the \emph{edge stalk} of $\mathcal{F}$ at $e$.
\item For each incidence of a vertex $v$ with an edge $e$, a \emph{restriction morphism} $\mathcal{F}_{v\face e}: \mathcal{F}(v) \to \mathcal{F}(e)$ in $\mathbf{C}$.
\end{itemize}
\end{definition}

We assume that $X$ has no loop edges, subdividing if necessary, so that each incidence $v \face e$ names a single endpoint of $e$. In our applications, $\mathbf{C}$ will typically be the category of groups ($\mathbf{Grp}$), abelian groups ($\mathbf{Ab}$), or sets ($\mathbf{Set}$). Stalks encode local data objects, while restriction morphisms encode how data propagates from vertices to edges.

Of particular importance is the \emph{constant sheaf} $\underline{G}$ associated with a group $G$. This assigns $\underline{G}(v) = \underline{G}(e) = G$ for all vertices $v$ and edges $e$, with all restriction maps being the identity homomorphism.

{\em Remark:} The base space $X$ represents the abstract topological spine of the configuration being analyzed. For instance, in the images of Figure \ref{fig:classics}, $X$ can be taken to be the circle $S^1$, represented as a graph with three or four vertices and edges, the vertices representing corners. This abstract space must be distinguished from its visual representation, such as a specific line drawing projected onto a 2-D plane.

{\em Notational Convention.} Throughout this paper, we represent the cyclic group of order two multiplicatively as $\mathbb{Z}_2 = \{+1, -1\}$ rather than additively. 

%------------------------------------------
\subsection{Network Sheaf Cohomology: Abelian}
%------------------------------------------

Sheaves -- local data structures -- are studied through their cohomology -- global invariants. The network case provides a discrete analogue of classical sheaf cohomology.

%\begin{definition}
Given a network sheaf $\mathcal{F}$ of abelian groups on a graph $X = (V, E)$, we define the cochain complex $C^\bullet(X; \mathcal{F})$ as:
\[
C^0(X; \mathcal{F}) = \prod_{v \in V} \mathcal{F}(v) 
\quad : \quad
C^1(X; \mathcal{F}) = \prod_{e \in E} \mathcal{F}(e)
\]

The coboundary map $\delta: C^0(X; \mathcal{F}) \to C^1(X; \mathcal{F})$ is defined as follows. For an oriented edge $e=(\vertexneg, \vertexpos)$, we set:
\begin{equation}
\label{eq:coboundary}
(\delta \xi)_e = \mathcal{F}_{\vertexpos\face e}(\xi_{\vertexpos}) - \mathcal{F}_{\vertexneg\face e}(\xi_{\vertexneg})
\end{equation}

The \emph{network sheaf cohomology} of $X$ with coefficients in $\mathcal{F}$ is defined as:
\begin{align*}
H^0(X; \mathcal{F}) &= \ker\,\delta \\
H^1(X; \mathcal{F}) &= \coker\,\delta = C^1(X; \mathcal{F})/\text{im}\,\delta
\end{align*}
%\end{definition}

These cohomology groups have clear interpretations:

\begin{itemize}
\item $H^0(X; \mathcal{F})$ consists of \emph{global sections} of $\mathcal{F}$ -- assignments of data to all vertices that can be consistently extended to all edges according to the sheaf structure. A global section represents a globally consistent state.

\item $H^1(X; \mathcal{F})$ measures the obstruction to extending locally compatible data to a global section. Non-trivial classes in $H^1$ characterize obstructions that prevent such a global extension.
\end{itemize}

For a constant sheaf $\underline{A}$ of an abelian group $A$, the coboundary map simplifies to:
\[
(\delta \xi)_e = \xi_{\vertexpos} - \xi_{\vertexneg}
\]
This corresponds to the usual simplicial or cellular coboundary operator.

%------------------------------------------
\subsection{Nonabelian Sheaf Cohomology}
\label{sec:nonabelian}
%------------------------------------------

When the stalks of $\mathcal{F}$ are nonabelian groups there is no coboundary homomorphism and no quotient group. First cohomology survives as a pointed set of orbits, and it classifies $\mathcal{F}$-torsors over $X$ (Section~\ref{sec:torsors}).

A 1-cochain $\eta \in C^1(X; \mathcal{F})$ assigns to each oriented edge $e$ an element $\eta_e \in \mathcal{F}(e)$. A graph has no 2-cells, hence no cocycle condition: every 1-cochain is a 1-cocycle, $Z^1(X; \mathcal{F}) = C^1(X; \mathcal{F})$. Cocycle conditions, and the obstructions in $H^2$ that they support, enter with 2-cells \cite{GhristGhrist2026Degrees}. The group $C^0(X;\mathcal{F})=\prod_{v}\mathcal{F}(v)$ acts on $C^1(X;\mathcal{F})$ by \emph{gauge transformations}: two 1-cocycles $\eta, \eta'$ are \emph{cohomologous} if there exists $\xi \in C^0(X; \mathcal{F})$ such that for each oriented edge $e = (\vertexneg, \vertexpos)$
\begin{equation}
\label{eq:gauge}
\eta'_e = \mathcal{F}_{\vertexneg\face e}(\xi_{\vertexneg})^{-1} \cdot \eta_e \cdot \mathcal{F}_{\vertexpos\face e}(\xi_{\vertexpos}) .
\end{equation}
The first cohomology $H^1(X; \mathcal{F}) = Z^1(X; \mathcal{F})/{\sim}$ is the set of gauge orbits, pointed by the orbit of the identity cochain: a pointed set, not a group. When the stalks are abelian and written additively, \eqref{eq:gauge} reads $\eta'_e = \eta_e + (\delta\xi)_e$, and the definition reduces to that of the previous subsection.

For a constant sheaf $\underline{G}$ the gauge action is
\[ \eta'_e = \xi_{\vertexneg}^{-1} \cdot \eta_e \cdot \xi_{\vertexpos} , \]
and $\eta$ lies in the distinguished orbit precisely when $\eta_e = \xi_{\vertexneg}^{-1} \xi_{\vertexpos}$ for some $\xi \in C^0(X;\underline{G})$. We write $(\delta \xi)_e = \xi_{\vertexneg}^{-1} \xi_{\vertexpos}$ for such a \emph{coboundary}, with the caution that $\delta$ is now a map of sets and not a homomorphism.

{\em Holonomy.} Paths are read in the order of traversal. An edge path $\gamma = e_1 e_2 \cdots e_m$ has holonomy
\[
\mathrm{hol}_\eta(\gamma) = \eta_{e_1}^{\pm 1}\,\eta_{e_2}^{\pm 1}\cdots\eta_{e_m}^{\pm 1},
\]
with exponent $-1$ on an edge traversed against its orientation. Concatenation of paths becomes multiplication, $\mathrm{hol}_\eta(\gamma\gamma') = \mathrm{hol}_\eta(\gamma)\,\mathrm{hol}_\eta(\gamma')$, and under a gauge transformation the product telescopes: a path from $x$ to $y$ has its holonomy replaced by $\xi_x^{-1}\,\mathrm{hol}_\eta(\gamma)\,\xi_y$, so the holonomy of a loop based at $x_0$ is conjugated by $\xi_{x_0}$.

\begin{theorem}[Classification via Fundamental Group]
\label{thm:h1_pi1}
For a graph $X$ representing a connected space with basepoint $x_0$ and a constant sheaf $\underline{G}$ of a group $G$, there is a natural bijection of pointed sets:
\[
H^1(X; \underline{G}) \cong \mathrm{Hom}(\pi_1(X,x_0), G)/G
\]
where the right side denotes conjugacy classes of homomorphisms $\pi_1(X,x_0) \to G$, with $G$ acting by conjugation. This bijection preserves the distinguished points: the trivial cohomology class corresponds to the conjugacy class of the trivial homomorphism.
\end{theorem}

{\em Proof sketch:} The bijection sends $[\eta]$ to the conjugacy class of $\gamma \mapsto \mathrm{hol}_\eta(\gamma)$, well defined by the telescoping above. Choose a spanning tree $T \subset X$. A gauge transformation, built outward from $x_0$ along $T$, sets $\eta_e = 1$ on every edge of $T$. The remaining edges correspond to free generators of $\pi_1(X, x_0)$, and their labels define an arbitrary homomorphism from this free group to $G$. The gauge freedom left after fixing $\eta|_T = 1$ is a single constant $\xi$, acting by simultaneous conjugation. \qed

This result is especially useful for analyzing visual paradoxes, as it allows us to characterize the obstruction directly in terms of representations of the fundamental group of the underlying space.

\begin{remark}
Two different choices hide in Theorem \ref{thm:h1_pi1}. Conjugation by $G$ accounts for the choice of origin in the fiber over the basepoint $x_0$: it is the gauge freedom that survives once a spanning tree has been fixed. Moving the basepoint itself, from $x_0$ to $x_1$, is a separate choice, implemented by transport along a path from $x_0$ to $x_1$; for connected $X$ the two descriptions agree, again up to conjugation.
\end{remark}

\begin{remark}
When $G$ is abelian, the conjugation action is trivial, and we have:
\[
H^1(X; \underline{G}) \cong \mathrm{Hom}(\pi_1(X,x_0), G) \cong \mathrm{Hom}(H_1(X), G)
\]
since every homomorphism to an abelian group factors through the abelianization $H_1(X)$ of $\pi_1(X,x_0)$.
\end{remark}

This framework provides the foundation for understanding how local geometric information can be consistently defined yet globally incompatible, creating the precise mathematical structure of visual paradoxes.

%%%%%%%%%%%%%%%%%%%%%%%%%%%%%%%%%%%%%%%%%%%%%
\section{Torsors \& Relative Change}
\label{sec:torsors}
%%%%%%%%%%%%%%%%%%%%%%%%%%%%%%%%%%%%%%%%%%%%%

The visual paradoxes examined in Section \ref{sec:gallery} each exhibit a tension between local consistency and global impossibility. This tension is naturally modeled through the mathematical concept of a $G$-torsor for some appropriate structure group $G$. In this section, we develop the mathematical foundations of this approach, introducing the concept of a network torsor as a discrete analogue of the classical torsor notion. There is nothing in this section that will be new to experts; for the typical reader, this gives a gentle self-contained introduction. 

%------------------------------------------
\subsection{Torsors: Intuition and Definition}
%------------------------------------------

Intuitively, a torsor is a space of states in which relative positions are well defined but no position is distinguished as the origin. Mathematically, a $G$-torsor over a space $X$ is a principal fiber bundle over $X$ with fiber $G$, where $G$ acts freely and transitively on each fiber. Classically, torsors are often introduced as {\em principal homogeneous spaces} -- spaces equipped with a free and transitive group action but without a canonical choice of origin. The critical property of torsors is that non-trivial torsors are precisely those that do not admit a global section.

This concept captures with precision the essence of visual paradoxes: local consistency and global inconsistency coexisting. In the Penrose staircase, for example, each local segment represents a well-defined height change (a relative attribute), but no globally consistent assignment of absolute heights exists.

%------------------------------------------
\subsection{Network $G$-Torsors}
\label{sec:network-torsors}
%------------------------------------------

We develop a discrete analogue of torsors for our combinatorial setting of graphs and networks. For more general approaches to torsors, see the classic references \cite{giraud1971cohomologie,kashiwara1990sheaves,maclane1992sheaves}. Our formulation is particularly tailored to computation and application. 

\begin{definition}
\label{def:torsor}
Let $G$ be a group and $X=(V,E)$ an oriented graph. A \emph{network $G$-torsor} $\torsor$ is a network sheaf of sets on $X$ 
equipped with a left action of $G$ on each stalk, $(g,p) \mapsto g \cdot p$, such that:

\begin{itemize}
\item For every vertex $v \in V$ and every edge $e \in E$, the stalks $\torsor(v)$ and $\torsor(e)$ are non-empty.

\item The actions are free and transitive: for any $p, p'$ in a common stalk, there exists a unique $g \in G$ such that $g \cdot p = p'$.

\item The restriction maps $\torsor_{v\face e}: \torsor(v) \to \torsor(e)$ are $G$-equivariant: $\torsor_{v\face e}(g \cdot p) = g \cdot \torsor_{v\face e}(p)$ for any $p \in \torsor(v)$ and $g \in G$.
\end{itemize}
\end{definition}

Each $g \in G$ acts by a bijection of each stalk, and once a point $p \in \torsor(x)$ has been chosen, the orbit map $g \mapsto g \cdot p$ is a bijection $G \to \torsor(x)$. No such choice is canonical, and that is the point of the definition: any two states in a stalk differ by a unique group element, while no state is distinguished as the origin. Free and transitive actions encode {\em relative positioning without absolute reference}. We use left actions throughout, so that restriction maps can be built from right multiplications, which commute with them (Construction \ref{constr:torsor_from_cocycle}).

\begin{definition}
\label{def:iso-torsor}
Let $\torsor$ and $\torsor'$ be network $G$-torsors over the same graph $X$.
A \emph{morphism of torsors} $\Phi:\torsor\to\torsor'$ consists of maps on stalks 
\[
\Phi_x : \torsor(x)\longrightarrow\torsor'(x)\qquad(x\in V\cup E)
\]
such that
\begin{enumerate}
\item $\Phi$ is $G$–equivariant: $\Phi_x(g\cdot p)=g\cdot\Phi_x(p)$ for every $p\in\torsor$ and $g\in G$;
\item it respects the restriction maps: for every incidence $(v\face e)$
      \[
      \Phi_e\!\bigl(\,\torsor_{v\face e}(p)\bigr)
        =\torsor'_{v\face e}\!\bigl(\,\Phi_v(p)\bigr).
      \]
\end{enumerate}
It is a fact about torsors that all morphisms are stalk-wise bijections and indeed are isomorphisms -- they admit an inverse torsor-morphism.
\end{definition}

\begin{definition}
\label{def:trivial-torsor}
A network $G$–torsor $\torsor$ over $X$ is \emph{trivial} if it admits a global section.
Equivalently, $\torsor$ is
trivial exactly when $\torsor\cong\underline G$, the constant torsor.
\end{definition}

A global section is a consistent reference frame throughout the network. We write $\Gamma(X;\torsor)$ for the set of global sections. A trivial torsor over a connected graph has one global section for each element of $G$, and none is preferred: triviality means that a compatible frame exists, not that one is given.

%------------------------------------------
\subsection{Classification of Network Torsors}
\label{sec:classification}
%------------------------------------------

The connection between torsors and first cohomology is made precise by the following standard result \cite{giraud1971cohomologie, stacks-project}, adapted to the network setting. It and Proposition \ref{prop:network_global_section} are included for their explicit computational form, in terms of edge labels and gauge transformations.

\begin{theorem}[Classification of Network Torsors]
\label{thm:network_torsor_classification}
Let $X = (V, E)$ be a connected graph and $G$ a group. There is a canonical bijection:
\begin{equation}
H^1(X; \underline{G}) 
\quad \longleftrightarrow  \quad
\{\text{isomorphism classes of network $G$-torsors over $X$}\}
\end{equation}

This bijection sends the distinguished element (identity class) of $H^1(X; \underline{G})$ to the class of the trivial network $G$-torsor.
\end{theorem}

{\em Proof:} Given a 1-cocycle $\eta$, Construction \ref{constr:torsor_from_cocycle} below produces a torsor $\torsor_\eta$ whose stalks are copies of $G$ with left multiplication and whose restriction maps are right multiplications: $(\torsor_\eta)_{\vertexneg\face e}(g) = g$ and $(\torsor_\eta)_{\vertexpos\face e}(g) = g\,\eta_e^{-1}$ for $e = (\vertexneg,\vertexpos)$. These commute with the left action, so $\torsor_\eta$ is a network $G$-torsor.

Conversely, let $\torsor$ be a network $G$-torsor, and choose a point $p_x$ in every stalk, vertex and edge alike. For $e = (\vertexneg,\vertexpos)$ there are unique $a_e, b_e \in G$ with $\torsor_{\vertexneg\face e}(p_{\vertexneg}) = a_e\cdot p_e$ and $\torsor_{\vertexpos\face e}(p_{\vertexpos}) = b_e \cdot p_e$. Replacing $p_e$ by $a_e \cdot p_e$ normalizes the first of these; identifying each stalk with $G$ via $g \mapsto g\cdot p_x$ then carries $\torsor$ to $\torsor_\eta$ with $\eta_e = a_e b_e^{-1}$. Replacing each vertex point $p_v$ by $\xi_v \cdot p_v$ replaces $\eta_e$ by $\xi_{\vertexneg}\,\eta_e\,\xi_{\vertexpos}^{-1}$, the gauge transformation \eqref{eq:gauge} by $\xi^{-1}$, so the class $[\eta]$ depends only on $\torsor$.

It remains to see that $\torsor_\eta \cong \torsor_{\eta'}$ exactly when $\eta$ and $\eta'$ are cohomologous. A map $G \to G$ equivariant for left multiplication is right multiplication by some element, so a morphism $\Phi: \torsor_\eta \to \torsor_{\eta'}$ has the form $\Phi_x(g) = g\,c_x$. Compatibility with the restriction maps at the tail of $e$ forces $c_e = c_{\vertexneg}$, and at the head it forces $\eta_e^{-1} c_{\vertexneg} = c_{\vertexpos}\,\eta'^{-1}_e$, that is, $\eta'_e = c_{\vertexneg}^{-1}\,\eta_e\,c_{\vertexpos}$: this is \eqref{eq:gauge} with $\xi = c$. Every such $\Phi$ is invertible. The identity cochain yields the constant torsor $\underline{G}$, the trivial torsor. \qed

\begin{proposition}[Global Sections of Network Torsors]
\label{prop:network_global_section}
A network $G$-torsor $\torsor$ over a graph $X = (V, E)$ admits a global section if and only if $\torsor$ is isomorphic to the trivial network torsor, i.e., its class in $H^1(X; \underline{G})$ is the distinguished element.
\end{proposition}

{\em Proof:} If $\torsor$ is the trivial network torsor, it is isomorphic to the constant sheaf $\underline{G}$ with its standard $G$-action, which clearly admits the identity section.

Conversely, given a global section $\glosection \in \Gamma(X; \torsor)$, we construct an isomorphism $\varphi: \underline{G} \to \torsor$ as follows. For each vertex $v \in V$, define $\varphi_v(g) = g \cdot \glosection_v$ for all $g \in G$. For each edge $e = (\vertexneg, \vertexpos)$, the compatibility condition for $\glosection$ ensures that $\torsor_{\vertexneg \face e}(\glosection_{\vertexneg}) = \torsor_{\vertexpos \face e}(\glosection_{\vertexpos})$. Call this common value $p_e \in \torsor(e)$. Define $\varphi_e(g) = g \cdot p_e$ for all $g \in G$.

The $G$-equivariance of the restriction maps and the free and transitive nature of the $G$-action ensure that $\varphi$ is an isomorphism of network $G$-torsors. \qed

For abelian structure groups, we have an additional characterization:

\begin{proposition}[Classification of Network Torsors with Abelian Structure Groups]
\label{prop:network_abelian_classification}
For a connected graph $X = (V, E)$ and an abelian group $A$, there is a canonical isomorphism:
\[
H^1(X; \underline{A}) \cong \mathrm{Hom}(H_1(X), A)
\]
where $H_1(X)$ is the first homology group of $X$.
\end{proposition}

This follows from the universal coefficient theorem and abelianization of $\pi_1$. 

%------------------------------------------
\subsection{Constructing Torsors from Cocycles}
%------------------------------------------

The classification theorem provides a method to construct network torsors from cohomology classes. Given a cocycle $\eta \in Z^1(X; \underline{G})$, we construct the torsor $\torsor_\eta$ as follows:

\begin{construction}[Torsor from Cocycle]
\label{constr:torsor_from_cocycle}
Let $\eta \in Z^1(X; \underline{G})$ be a 1-cocycle. We define a network $G$-torsor $\torsor_\eta$ as follows:
\begin{itemize}
\item For each vertex $v \in V$, set $\torsor_\eta(v) = G$ as a set.
\item For each edge $e \in E$, set $\torsor_\eta(e) = G$ as a set.
\item The $G$-action on each stalk is by left multiplication: $g \cdot p = gp$.
\item For an edge $e = (\vertexneg, \vertexpos)$, define the restriction maps:
\[
(\torsor_\eta)_{\vertexneg \face e}(g) = g 
\quad : \quad
(\torsor_\eta)_{\vertexpos \face e}(g) = g \cdot \eta_e^{-1}
\]
\end{itemize}
\end{construction}

Right multiplication commutes with left multiplication, so the restriction maps are equivariant. A global section is an assignment $s_v \in G$ with $s_{\vertexneg} = s_{\vertexpos}\,\eta_e^{-1}$ on each edge, that is,
\[ s_{\vertexpos} = s_{\vertexneg}\,\eta_e . \]
The cocycle is forward transport, and along an edge path it composes in the order of traversal, as in the holonomy convention of Section \ref{sec:nonabelian}. For an abelian group written additively this is the familiar height rule $s_{\vertexpos} = s_{\vertexneg} + \eta_e$.

Reversing the orientation of an edge and inverting its label gives an isomorphic torsor; the class in $H^1$ is consequently independent of the orientations chosen.
This construction makes the correspondence between cohomology classes and torsors explicit. The cocycle $\eta$ encodes the ``twist'' or discrepancy that prevents the existence of a global section when $\eta$ is non-trivial.

Specifically, a cocycle $\eta \in Z^1(X; \underline{G})$ represents the trivial class in $H^1(X; \underline{G})$ if and only if it is cohomologous to the identity cocycle, which, according to the equivalence relation, means there exists a 0-cochain $\xi \in C^0(X; \underline{G})$ such that $\eta_e = \xi_{\vertexneg}^{-1} \cdot \xi_{\vertexpos}$ for every oriented edge $e=(\vertexneg, \vertexpos)$.

For a 1-cocycle to be trivial, it must be a coboundary: $\eta = \delta \xi$ for some 0-cochain $\xi \in C^0(X; \underline{G})$. In this case, a global section of $\torsor_\eta$ can be constructed as $s_v = \xi_v$ for each vertex $v$.

%------------------------------------------
\subsection{Interpretation Problems}
\label{sec:interpretation}
%------------------------------------------

A torsor is the right object when the relative transformations act on the local states freely and transitively. Often they do not: reflections act on heights with fixed points. The models of this paper are instances of the following structure.

\begin{definition}
\label{def:interpretation}
An \emph{interpretation problem} on a connected graph $X$ consists of a group $G$, a transport cocycle $\eta \in C^1(X;\underline{G})$, and a set $S$ of local states with a left action of $G$. A \emph{consistent interpretation} assigns a state $\sigma_v \in S$ to each vertex so that $\sigma_{\vertexneg} = \eta_e \cdot \sigma_{\vertexpos}$ on every edge $e = (\vertexneg, \vertexpos)$; equivalently, it is a global section of the sheaf of $G$-equivariant maps from $\torsor_\eta$ to $S$.
\end{definition}

\begin{lemma}
\label{lem:fixed-point}
A consistent interpretation exists if and only if the image of $\mathrm{hol}_\eta$ in $G$ has a common fixed point in $S$. When $G$ acts freely and transitively on $S$, this holds if and only if $[\eta]$ is the distinguished class of $H^1(X;\underline{G})$.
\end{lemma}

{\em Proof:} Along any edge path $\gamma$ from $x_0$ to $x$, the condition gives $\sigma_{x_0} = \mathrm{hol}_\eta(\gamma)\cdot\sigma_x$, so an interpretation is determined by $\sigma_{x_0}$, and it is well defined exactly when $\sigma_{x_0} = \mathrm{hol}_\eta(\gamma)\cdot\sigma_{x_0}$ for every loop $\gamma$ at $x_0$. For a free action a fixed point exists only for the trivial holonomy, which by Theorem \ref{thm:h1_pi1} is the distinguished class. \qed

When the states are coordinates on local charts, $\eta_e$ is the change of coordinates from the chart at the head of $e$ to the chart at its tail, and a consistent interpretation is a single state read in every chart; the relation $\sigma_{x_0} = \mathrm{hol}_\eta(\gamma)\cdot\sigma_x$ then says that $\mathrm{hol}_\eta(\gamma)$ carries coordinates at the end of $\gamma$ to coordinates at its start. For $S = G$ with left multiplication, $\sigma_v = s_v^{-1}$ recovers the sections of $\torsor_\eta$ in Construction \ref{constr:torsor_from_cocycle}.

When $S$ is a torsor, a nontrivial class forbids a consistent interpretation. When the action has stabilizers, it need not: a nontrivial class whose holonomy fixes a common state admits one. The M\"obius staircase of Section \ref{sec:non_orientable} is the example to keep in mind: for even $\Delta$, heights alone admit a consistent assignment, whereas height frames do not.

%%%%%%%%%%%%%%%%%%%%%%%%%%%%%%%%%%%%%%%%%%%%%
\section{Paradoxical Figures: A Cohomological Analysis}
\label{sec:examples}
%%%%%%%%%%%%%%%%%%%%%%%%%%%%%%%%%%%%%%%%%%%%%

We now apply the framework developed in Sections \ref{sec:toolkit}-\ref{sec:torsors} to analyze several classes of impossible figures. For each case, we identify the underlying topological space, the appropriate structure group, and the corresponding torsor that captures the essential visual paradox. Before examining specific examples, we establish a common analytical approach:

\begin{enumerate}
    \item \textbf{Network Structure}: Model the figure's local constraints as a graph, typically homotopy equivalent to a simple space ($S^1$, $S^1 \vee S^1$, etc.). When the figure lives on a surface, we write $\Lambda$ for the constraint graph and $X$ for the surface; data on $\Lambda$ define a class on $X$ only once they are shown to extend across its 2-cells.
    
    \item \textbf{Structure Group}: Identify the appropriate structure group $G$ representing the local ambiguity (height changes, orientation flips, etc.)
    
    \item \textbf{Constraint Cocycle}: Encode the perceived local transformations as a cocycle $\eta \in Z^1(X; \underline{G})$.
    
    \item \textbf{Torsor Classification}: By Construction \ref{constr:torsor_from_cocycle}, from $\eta$ determine the corresponding torsor $\torsor_\eta$ and its non-triviality via $H^1(X; \underline{G})$.
\end{enumerate}

Every model below is an interpretation problem in the sense of Definition \ref{def:interpretation}. It names a base, a structure group, a set of local states, and transport, prescribed and locally constant along edges, with coefficient groups treated as abstract groups. When the local states form a torsor, a paradox is a nontrivial class, and by Proposition \ref{prop:network_global_section} a consistent interpretation exists precisely when $[\eta]$ is the distinguished class. When the action has stabilizers, Lemma \ref{lem:fixed-point} is the criterion, and nontriviality alone forbids nothing. A visual paradox is the absence of a consistent interpretation for the specified problem.

%------------------------------------------
\subsection{Height-Based Paradoxes: The Penrose Staircase}
\label{sec:height_paradoxes}
%------------------------------------------

The Penrose staircase (Figure \ref{fig:classics2}) provides the framework for height-based paradoxes.

{\em Network Structure.} The Penrose staircase is modeled as a graph $X$ with vertices $V = \{v_1, v_2, v_3, v_4\}$ representing the corners and edges $E = \{e_1, e_2, e_3, e_4\}$ representing the flights of stairs, where $e_i$ connects $v_i$ and $v_{i+1}$ (with indices taken modulo 4). This graph $X$ is homotopy equivalent to the circle $S^1$.

{\em Structure Group.} The attribute in this paradox is discrete height change represented by the additive group $\mathbb{Z}$. We define the structure sheaf as the constant sheaf $\mathcal{G} = \underline{\mathbb{Z}}$ over $X$. This is a modeling choice -- we could alternatively use the real numbers $\mathbb{R}$ for height or the positive reals $\mathbb{R}^+$ for multiplicative height scaling.

{\em Constraint Cocycle.} Traversing each edge $e_i$ corresponds to a perceived height change $\Delta_i \in \mathbb{Z}$. In the Penrose staircase of Figure \ref{fig:classics2} these $\Delta_i$ are all positive (or all negative). These perceived changes are encoded by a 1-cochain $\eta \in C^1(X; \underline{\mathbb{Z}})$ by setting $\eta_{e_i} = \Delta_i$.

{\em Cohomological Analysis.} The cochain $\eta$ must be a coboundary for a globally consistent height assignment to exist. That is, there must exist a 0-cochain $\xi \in C^0(X; \underline{\mathbb{Z}})$ representing absolute height assignments $\xi_{v_i}$ at each vertex such that $\eta = \delta \xi$. This requires:
\[
\eta_{e_i} = (\delta \xi)_{e_i} = \xi_{v_{i+1}} - \xi_{v_i}
\]

Summing these constraints around the complete loop
$\sum_{i=1}^4 \eta_{e_i} = \sum_{i=1}^4 (\xi_{v_{i+1}} - \xi_{v_i})$ yields the constraint $\sum_{i=1}^4 \Delta_i = 0$. This contradicts the visual perception of the Penrose staircase where $k = \sum_{i=1}^4 \Delta_i \neq 0$.

{\em Torsor Interpretation.} Since $k \neq 0$ the cocycle $\eta$ is not a coboundary and the corresponding torsor $\torsor_\eta$ is non-trivial and lacks a global section. The non-triviality of this torsor {\it is} the mathematical expression of the visual paradox. The obstruction is classified by $k \in \mathbb{Z} \cong H^1(S^1; \mathbb{Z})$, which quantifies the ``degree'' of the paradox in terms of the net height change accumulated around the loop.

Our visual system attempts to integrate these local height cues into a global model, but the underlying mathematical obstruction $[\eta] \in H^1(X; \underline{\mathbb{Z}})$ prevents a consistent interpretation. The impossibility is not merely a failure of our perception but reflects a genuine topological obstruction.

%------------------------------------------
\subsection{Three-Dimensional Cubic Staircases}
\label{sec:cubic_staircases}
%------------------------------------------

%We begin with perhaps the most intuitive class of visual paradoxes: three-dimensional cubic staircases where the obstruction arises directly from spatial translations.

{\em Network Structure.} Consider a cycle of cubes arranged as in Figure \ref{fig:penrose-stairs} where each cube shares exactly two faces with neighbors. The underlying graph $X$ has vertices at centers of the cubes and edges for adjacency relations, forming a cycle homotopy equivalent to $S^1$.

{\em Structure Group.} Since each cube-to-cube transition represents a unit translation in one of three coordinate directions, the structure group is $G = \mathbb{Z}^3$. The structure sheaf is constant $\underline{\mathbb{Z}^3}$ over $X$.

{\em Constraint Cocycle.} Each edge $e_i$ in the cycle corresponds to a unit translation $\vect{v}_i$, one of $\pm\vect{e}_1, \pm\vect{e}_2, \pm\vect{e}_3$, read from the drawing. A consistent arrangement of cubes would require the sum of these translations to vanish. The local translations are encoded by a 1-cochain $\eta \in C^1(X; \underline{\mathbb{Z}^3})$ with $\eta_{e_i} = \vect{v}_i$. The holonomy around the cycle is
$\vect{h} = \sum \vect{v}_i$.

{\em Cohomological Analysis.} For a globally consistent cube arrangement to exist, $\eta$ must be a coboundary with $\eta = \delta \xi$ for a $0$-cochain of absolute positions $\xi \in C^0(X; \underline{\mathbb{Z}^3})$. It is not, and the drawings show why.

{\em Why the Drawing Closes.} Each figure is the image of a chain of cubes under a linear projection $P : \mathbb{R}^3 \to \mathbb{R}^2$, and the drawn chain closes up in the plane. Hence $P(\vect{h}) = \sum P(\vect{v}_i) = 0$: the holonomy lies in the lattice $\ker P \cap \mathbb{Z}^3$. For the isometric projection of the top row of Figure \ref{fig:penrose-stairs} this lattice is generated by $(1,1,1)$; for the dimetric projection of the bottom row it is generated by $(2,2,1)$, since there a unit vertical step has the same image as two unit steps along one horizontal axis followed by two backward along the other. Each holonomy is therefore an integer multiple $k$ of one primitive vector, and only the multiple remains to be counted. \emph{The paradox hides in the kernel of the camera.} Counting unit steps around each loop gives\footnote{Figures are read left to right, then top to bottom. Each loop is traversed counterclockwise in the image, and in each row the basis of unit steps is oriented so that the primitive kernel vector has positive entries.}
\begin{equation*}
\label{eq:cubical-holos}
    \begin{pmatrix}2\\2\\2\end{pmatrix}
    \quad : \quad
    \begin{pmatrix}1\\1\\1\end{pmatrix}
    \quad : \quad
    \begin{pmatrix}1\\1\\1\end{pmatrix}
    \quad : \quad
    \begin{pmatrix}4\\4\\4\end{pmatrix}
    \quad : \quad
    \begin{pmatrix}-2\\-2\\-1\end{pmatrix}
    \quad : \quad
    \begin{pmatrix}-2\\-2\\-1\end{pmatrix}
    \quad : \quad
    \begin{pmatrix}2\\2\\1\end{pmatrix}
    \quad : \quad
    \begin{pmatrix}-2\\-2\\-1\end{pmatrix}
\end{equation*}

Since the corresponding cocycles are not coboundaries, the associated $\mathbb{Z}^3$-torsors are non-trivial.

{\em Torsor Interpretation.} The holonomy vector $\vect{h}$ measures how the structure fails to close properly in 3-D space following adjacency rules. The torsors, classified by elements of $H^1(S^1; \underline{\mathbb{Z}^3}) \cong \mathbb{Z}^3$, mathematically express the paradoxes.

We note that one could project these $\mathbb{Z}^3$-torsors onto a single dimension (such as height) to obtain simpler $\mathbb{Z}$-torsors. This dimensional reduction loses information about the spatial translation. %These different holonomy vectors indicate qualitatively distinct paradoxes, even when projected onto the same single dimension torsor.

%This three-dimensional perspective provides a natural foundation for understanding more specialized paradoxes focused on particular aspects of spatial perception, such as the height-based examples we consider next.

%------------------------------------------
\subsection{Identification: Cylindrical Staircase}
\label{sec:cylindrical_staircase}
%------------------------------------------

% The cylindrical staircase is a variation to the Penrose staircase that reflects how paradoxes can arise in different topological settings while maintaining the same structure.

{\em Network Structure.} The cylindrical staircase (Figure \ref{fig:penrose_cyl}[right]) is modeled as a graph $X$ with vertices at key positions along the cylinder and edges representing the connecting stair segments. Topologically, $X$ is again homotopy equivalent to $S^1$.

{\em Identification Space and Coordinates.} The cylinder is realized as an identification space. The cylinder is represented as a rectangle with its left and right edges identified, preserving orientation and height. 

Concretely, we may establish a local coordinate system $(x,y,z)$ where $x$ is the coordinate along the circumference of the cylinder, $y$ is the radial coordinate, and $z$ represents height, then the identification operates only on the $x$-coordinate. For points $(0,y,z)$ on the left edge and $(L,y,z)$ on the right edge (where $L$ is the width of the rectangle), the identification equates $(0,y,z) \sim (L,y,z)$ for all $y$ and $z$. Crucially, the $z$-coordinate (height) is preserved under this identification.

{\em Structure Group and Cocycle.} As with the Penrose staircase, we use the structure group $G = \mathbb{Z}$ to represent height changes, with structure sheaf $\mathcal{G} = \underline{\mathbb{Z}}$ with a 1-cochain of height changes $\eta \in C^1(X; \underline{\mathbb{Z}})$.

Paths that remain within the fundamental domain (the rectangle before identification) encounter no contradictions. It is only when a path crosses the identification boundary that the paradox emerges. For such a path, the net height change $k = \sum_i \eta_{e_i} \neq 0$ contradicts the constraint that a globally consistent height assignment would require $k = 0$.

{\em Torsor Classification.} The mathematical classification of this paradox is identical to that of the Penrose staircase: a non-trivial $\mathbb{Z}$-torsor over $S^1$ with obstruction class $k \in H^1(S^1; \mathbb{Z})$. This equivalence demonstrates that the essential mathematical structure of the paradox is independent of its specific visual realization.

%------------------------------------------
\subsection{Multiple Loops: Torus Staircases}
\label{sec:torus_paradoxes}
%------------------------------------------

The torus provides a richer setting with multiple independent loops, leading to more complex paradoxical structures (Figure \ref{fig:torus}).

{\em Network Structure.} Here the constraint graph and the ambient surface must be kept apart. In its simplest manifestation, a staircase on a torus is modeled by a constraint graph $\Lambda = S^1 \vee S^1$ carrying the two generating loops $a$ and $b$.\footnote{This holds for Figure \ref{fig:torus}[left]. For Figure \ref{fig:torus}[right] the constraint graph is larger, and its height data reduce to two torus periods only after one checks that the holonomy vanishes around every cycle of the graph that bounds in the torus.} The fundamental group of the graph is free, $\pi_1(\Lambda) \cong F_2 = \langle a, b \rangle$. The torus $T^2$ is obtained by attaching a 2-cell along the word $aba^{-1}b^{-1}$, so that $\pi_1(T^2) \cong \mathbb{Z}\times\mathbb{Z}$, and data on $\Lambda$ extend across this 2-cell precisely when the holonomy around the attaching word is trivial: $\mathrm{hol}(a)$ and $\mathrm{hol}(b)$ must commute.

{\em Torsor Analysis.} With the structure group $\mathbb{Z}$ for height, the commutation condition holds automatically, and
\[
H^1(\Lambda; \mathbb{Z}) \cong \mathrm{Hom}(F_2, \mathbb{Z}) \cong \mathbb{Z} \times \mathbb{Z} \cong \mathrm{Hom}(\mathbb{Z} \times \mathbb{Z}, \mathbb{Z}) \cong H^1(T^2; \mathbb{Z}) .
\]
The graph and the torus share their first cohomology with abelian coefficients because every homomorphism to an abelian group kills the commutator -- not because their fundamental groups agree, which they do not.

Each pair $(\Delta_a, \Delta_b) \in \mathbb{Z} \times \mathbb{Z}$ of height changes along the two fundamental loops corresponds to a distinct class in $H^1(X; \mathbb{Z})$ and classifies a possibly paradoxical figure.\footnote{Reversing the orientation of both loops replaces $(\Delta_a, \Delta_b)$ by $(-\Delta_a, -\Delta_b)$. This is pullback along a homeomorphism of the base -- a change of model -- and not an isomorphism of torsors over a fixed base.} The torus example introduces an important feature: multiple independent paths between points can result in different height changes, an impossibility in physically realizable structures.

%------------------------------------------
\subsection{Orientation Paradoxes: Alternating Depth}
\label{sec:binary_paradoxes}
%------------------------------------------

Not all paradoxes involve continuous quantities like height. The zigzag depth paradox (Figure \ref{fig:zigzag}) introduces binary perceptual ambiguity.

{\em Structure and Analysis.} For a zigzag path with $n$ corners arranged in a loop, we use the structure group $\mathbb{Z}_2 = \{+1, -1\}$ to model the binary depth ambiguity at each corner. The constraint cocycle $\eta$ encodes depth flips with $\eta_{e_i} = -1$ for all edges. For a global depth assignment to exist, we need a 0-cochain $\xi$ with $\xi_{v_{i+1}} = -\xi_{v_i}$ for each $i$. Following these constraints around the loop gives $\xi_{v_1} = (-1)^n \cdot \xi_{v_1}$, which is consistent if and only if $n$ is even.

For a zigzag with an odd number of corners, the corresponding $\mathbb{Z}_2$-torsor is non-trivial, classified by the generator of $H^1(S^1; \mathbb{Z}_2) \cong \mathbb{Z}_2$. This captures the impossibility of assigning consistent depth interpretations while maintaining local coherence around the loop.

%------------------------------------------
\subsection{Non-orientable Surfaces: M\"obius and Projective Staircases}
\label{sec:non_orientable}
%------------------------------------------

The M\"obius staircase (Figure \ref{fig:mobius}) and projective plane staircase (Figure \ref{fig:RP^2}) introduce non-orientability.

{\em Structure Group.} We use the semi-direct product $G = \mathbb{Z} \rtimes \mathbb{Z}_2$, where $\mathbb{Z}$ represents height and $\mathbb{Z}_2 = \{+1, -1\}$ represents orientation. Flipping orientation reverses the perceived direction of height change: climbing on the top surface becomes descending when viewed from the bottom. The group operation is $(h_1, \epsilon_1) \cdot (h_2, \epsilon_2) = (h_1 + \epsilon_1 h_2, \epsilon_1 \epsilon_2)$, and $G \cong D_\infty$ is the group of affine symmetries of the integer lattice of tread heights.

{\em Local States.} The group acts on heights by $(h,\epsilon) \cdot z = h + \epsilon z$. In the language of Definition \ref{def:interpretation}, the label $\eta_e$ acts as the change of coordinates from the height chart at the head of $e$ to the chart at its tail, so the translation part of $\mathrm{hol}_\eta(\gamma)$ is the height of the chart at the end of $\gamma$, read in the chart at its start; this is the sense of every height component quoted below. The action is not free: the reflection $(h,-1)$ fixes the height $h/2$, which is a tread height when $h$ is even. The torsor is the set of \emph{height frames} -- an origin on the lattice of treads together with a sense of up -- on which $G$ acts simply transitively. Heights alone form an associated space, not a torsor.

{\em Constraint Cocycle.} For both examples, modeling the core path as a triangle graph $X = C_3$ with a flat segment, a stair segment (height change $\Delta \neq 0$), and a twist segment (orientation flip), the cocycle $\eta$ yields values:
\begin{itemize}
    \item $\eta_{e_1} = (0, +1)$: Flat segment
    \item $\eta_{e_2} = (\Delta, +1)$: Stair segment
    \item $\eta_{e_3} = (0, -1)$: Twist segment with orientation flip
\end{itemize}

{\em Holonomy Analysis.} The holonomy around the loop is:
\[
g = \eta_{e_1} \cdot \eta_{e_2} \cdot \eta_{e_3} = (0, +1) \cdot (\Delta, +1) \cdot (0, -1) = (\Delta, -1)
\]

This holonomy is non-trivial, classifying a non-trivial torsor. Walking twice around the M\"obius staircase returns to the original position, height, and orientation, evident by $g^2 = (\Delta, -1) \cdot (\Delta, -1) = (0, +1)$, the identity element. For the projective plane, whose 2-cell is attached along the square of the core loop, this identity is precisely the condition for the model to extend across the 2-cell.

Two cautions follow. First, the height component of a reflection is not invariant under change of gauge: conjugation by the translation $(k,+1)$ gives
\[
(k,+1)(\Delta,-1)(k,+1)^{-1} = (\Delta + 2k, -1).
\]
The value $(\Delta,-1)$ is the holonomy based at the foot of the flight; a walker starting at step $j$ reads $(\Delta - 2j, -1)$. The class in $H^1(S^1;\underline{G})$ records the flip, the torsion, and $\Delta \bmod 2$, and nothing more. Second, heights alone see only the parity of $\Delta$. The reflection $z \mapsto \Delta - z$ fixes $\Delta/2$, so a consistent assignment of signed tread heights exists precisely when $\Delta$ is even; for odd $\Delta$ the height-only model is already obstructed. The height-frame model is obstructed in either case: what fails is a consistent height frame, and in particular a consistent sense of up (Remark \ref{rem:orientation-preserving}).

Despite their different ambient topologies, the M\"obius and projective plane examples exhibit the same paradoxical behavior when restricted to their representative loops: both are classified by reflection classes in $H^1(S^1; \underline{\mathbb{Z} \rtimes \mathbb{Z}_2})$.

%------------------------------------------
\subsection{Nonabelian Paradox: The Klein Ladder}
\label{sec:nonabelian_paradox}
%------------------------------------------

The Klein ladder (Figure \ref{fig:klein}) is a nonabelian visual paradox.

{\em Structure and Homomorphism.} Model the ladder by a constraint graph $\Lambda = S^1 \vee S^1$ carrying the loops $a$ and $b$, inside the Klein bottle $K$ obtained by attaching a 2-cell along $aba^{-1}b$. With the structure group $G = \mathbb{Z} \rtimes \mathbb{Z}_2 \cong D_\infty$ of the M\"obius example, the net transformations along the two loops, read from a chosen basepoint and counted in ladders, define a homomorphism $\rho: \pi_1(\Lambda) = F_2 \to G$:
\begin{itemize}
    \item $\rho(a) = (\alpha, -1)$ (height change $\alpha$, orientation flip)
    \item $\rho(b) = (\beta, +1)$ (height change $\beta \neq 0$, no orientation flip)
\end{itemize}

Any assignment on the free generators is a homomorphism of $F_2$. The model extends across the 2-cell of $K$ precisely when $\rho$ kills the attaching word, that is, when $\rho$ factors through $\pi_1(K) = \langle a, b \mid ab = b^{-1}a \rangle$. Using $(h, \epsilon)^{-1} = (-\epsilon h, \epsilon)$,
\begin{align*}
\rho(a)\rho(b) &= (\alpha, -1) \cdot (\beta, +1) = (\alpha - \beta, -1) \\
\rho(b)^{-1}\rho(a) &= (-\beta, +1) \cdot (\alpha, -1) = (\alpha - \beta, -1) ,
\end{align*}
so the relation holds for every $\alpha$ and $\beta$: a flip followed by a climb is a climb reversed followed by the flip. Note also that $\rho(a)^2 = (\alpha - \alpha, +1) = 1$. The representation therefore factors further, through the quotient $\pi_1(K) \twoheadrightarrow D_\infty$ that imposes $a^2 = 1$, and its image $\{(\alpha + m\beta, -1), (m\beta, +1) : m \in \mathbb{Z}\}$ is itself infinite dihedral. The ladder records a representation of the Klein-bottle group, not a faithful copy of it. From the basepoint marked in Figure \ref{fig:klein}, $\alpha = 2$ and $\beta = 1$.

{\em Invariants.} Conjugation by the translation $(k,+1)$ replaces $\alpha$ by $\alpha + 2k$ and fixes $\beta$; conjugation by the reflection $(0,-1)$ changes the signs of both. The class $[\rho]$ is therefore determined by $|\beta|$ and $\alpha \bmod 2$. As in the M\"obius example, the height change along the flipping loop $a$ depends on where one starts; the height change along $b$, a translation, does not, though its sign reverses with the sense of up.

{\em Nonabelian Nature.} Compare the two orders of traversal, read from the same basepoint:
\[
\rho(ab) = \rho(a)\rho(b) = (\alpha - \beta, -1), \qquad \rho(ba) = \rho(b)\rho(a) = (\alpha + \beta, -1).
\]
Both paths end flipped, at heights differing by $2\beta \neq 0$, so $\rho(ab) \neq \rho(ba)$. The contrast between the relator and the commutator makes the point visible:
\[
\rho(aba^{-1}b) = 1, \qquad \rho(aba^{-1}b^{-1}) = (-2\beta, +1) \neq 1 .
\]
The first loop closes up, as it must on the Klein bottle; the second, which would close up on a torus, returns displaced in height by $2\beta$. Conjugation preserves commutativity, so no change of gauge makes the image of $\rho$ abelian: in this sense the paradox is nonabelian. This does not mean that no abelian invariant detects it. The sign homomorphism $D_\infty \to \mathbb{Z}_2$ already carries $[\rho]$ to a nontrivial class. The nonabelian content lies in the full transport, which no single loop can exhibit, since the holonomy along one circle generates a cyclic group.

%%%%%%%%%%%%%%%%%%%%%%%%%%%%%%%%%%%%%%%%%%%%%
\section{Paradoxes \& Boundary Conditions}
\label{sec:boundary}
%%%%%%%%%%%%%%%%%%%%%%%%%%%%%%%%%%%%%%%%%%%%%

Loops are not the only source of paradox. A path graph $P_n$ is contractible, so $H^1(P_n;\underline{G})$ is trivial for every group $G$: every torsor with constant structure group over a path has a global section. Yet the Necker sequence of Figure \ref{fig:necker_gradient} and the impossible bars of Figure \ref{fig:boundary} are paradoxical, because their ends are pinned by visual cues. The obstruction is relative to the pinned set, and it can be written in two equivalent ways: as relative cohomology, and as a torsor for a structure sheaf that is trivial on the pinned set.

%------------------------------------------
\subsection{Relative Cohomology}
%------------------------------------------

Let $X$ be a graph and $A \subset X$ a subcomplex; in our examples $A$ is a set of pinned vertices. For an abelian group $B$, the relative cochains $C^\bullet(X, A; B)$ are the cochains that vanish on the cells of $A$, with the coboundary \eqref{eq:coboundary}. The long exact sequence of the pair begins
\begin{equation}
\label{eq:LES}
H^0(X; B) \xrightarrow{\res} H^0(A; B) \xrightarrow{\delta^*} H^1(X, A; B) \to H^1(X; B) .
\end{equation}
Here $\res$ restricts global sections, and the connecting map $\delta^*$ sends boundary data $\beta \in H^0(A;B)$ to the class of $\delta\tilde\beta$, where $\tilde\beta$ extends $\beta$ by zero off $A$. By exactness, $\beta$ extends to a global section precisely when $\delta^*\beta = 0$.

For a path $X = P_n$ with vertices $v_0, \ldots, v_n$ and pinned ends $A = \{v_0, v_n\}$, summing a relative 1-cocycle along the path identifies $H^1(X,A;B) \cong B$, and $\delta^*(b_0, b_n) = b_n - b_0$: the connecting class of the boundary data is their difference. For the Necker sequence, $B = \mathbb{Z}_2$, written multiplicatively, and the contradictory cues $\beta = (+1,-1)$ give the nontrivial class.

%------------------------------------------
\subsection{Inhomogeneous Network Torsors}
%------------------------------------------

To capture boundary constraints within the torsor framework, we replace the constant sheaf $\underline{G}$ with a sheaf of groups $\mathcal{G}$ that varies across the network.

\begin{definition}
Let $X=(V,E)$ be an oriented graph and $\mathcal{G}$ a network sheaf of groups on $X$. A \emph{network $\mathcal{G}$-torsor} $\torsor$ consists of:
\begin{enumerate}
\item A network sheaf of non-empty sets $\torsor$ on $X$.
\item A family of left group actions $\mathcal{G}(v) \times \torsor(v) \to \torsor(v)$ and $\mathcal{G}(e) \times \torsor(e) \to \torsor(e)$, denoted $(g, p) \mapsto g \cdot p$.
\item These actions are free and transitive: for any $p, p' \in \torsor(v)$, there exists a unique $g \in \mathcal{G}(v)$ such that $g \cdot p = p'$, and similarly for edges.
\item The restriction maps $\torsor_{v\face e}: \torsor(v) \to \torsor(e)$ are compatible with the group actions via the sheaf restriction maps $\mathcal{G}_{v \face e}: \mathcal{G}(v) \to \mathcal{G}(e)$:
\[
\torsor_{v\face e}(g \cdot p) = \mathcal{G}_{v \face e}(g) \cdot \torsor_{v\face e}(p)
\]
\end{enumerate}
\end{definition}

\begin{theorem}[Classification of Network Torsors over Structure Sheaves]
\label{thm:nonconstant_classification}
Let $X = (V, E)$ be a graph and $\mathcal{G}$ a sheaf of groups on $X$. There is a canonical bijection:
\begin{equation}
H^1(X; \mathcal{G}) \longleftrightarrow \{\text{isomorphism classes of network $\mathcal{G}$-torsors over $X$}\}
\end{equation}
This bijection sends the distinguished element of $H^1(X; \mathcal{G})$ to the class of the trivial network $\mathcal{G}$-torsor.
\end{theorem}

The proof is that of Theorem \ref{thm:network_torsor_classification}, with the maps $\mathcal{G}_{v\face e}$ inserted where needed. The torsor of a cocycle $\eta \in C^1(X;\mathcal{G})$ has stalks $\mathcal{G}(x)$ with left multiplication and restriction maps $(\torsor_\eta)_{\vertexneg\face e}(g) = \mathcal{G}_{\vertexneg\face e}(g)$ and $(\torsor_\eta)_{\vertexpos\face e}(g) = \mathcal{G}_{\vertexpos\face e}(g)\,\eta_e^{-1}$, and cocycles give isomorphic torsors exactly when they are related by \eqref{eq:gauge}.

%------------------------------------------
\subsection{Structure Sheaves for Boundary Constraints}
%------------------------------------------

\begin{definition}
Let $G$ be a group and $A \subset X$ a subcomplex. The \emph{pinned sheaf} $\mathcal{G}_A$ is the subsheaf of $\underline{G}$ with stalk $\{1\}$ on the cells of $A$ and $G$ on all other cells; its restriction maps are the identity between cells outside $A$ and the inclusion of the identity element on cells of $A$.
\end{definition}

For abelian coefficients the pinned sheaf and relative cohomology are the same thing, already at the level of cochains.

\begin{proposition}
\label{prop:pinned-relative}
For an abelian group $B$, $C^\bullet(X; B_A) = C^\bullet(X, A; B)$ as cochain complexes. Consequently $H^1(X; B_A) \cong H^1(X, A; B)$.
\end{proposition}

{\em Proof:} The stalks of $B_A$ vanish exactly on the cells of $A$, so its cochains are the cochains vanishing on $A$. A restriction map out of a cell of $A$ is the zero map, so the coboundary \eqref{eq:coboundary} of $B_A$ is the relative coboundary. \qed

For a nonabelian group the right-hand side has no ready-made meaning, and we take the left-hand side as its definition: $H^1(X,A;\underline{G}) := H^1(X;\mathcal{G}_A)$, the set of 1-cocycles of $\mathcal{G}_A$, which equal $1$ on the edges of $A$, modulo gauge transformations \eqref{eq:gauge} by 0-cochains $\xi$ with $\xi = 1$ on $A$. This is the relative gauge model, and it is pointed by the class of the identity cochain.

{\em Pinned Torsors.} Boundary cues enter as a torsor. Let $A$ be a set of vertices. Fix a prescribed transport $\eta \in C^1(X;\underline{G})$ and boundary data $\beta = (\beta_a)_{a \in A}$ with $\beta_a \in G$. The \emph{pinned torsor} $\torsor_{\eta,\beta}$ has stalk $G$ on every cell outside $A$ and the singleton $\{\beta_a\}$ at each $a \in A$. Its restriction maps are those of $\torsor_\eta$ (Construction \ref{constr:torsor_from_cocycle}), applied to $\beta_a$ at a pinned vertex: the singleton maps to $\beta_a$ in the stalk of an edge with tail $a$, and to $\beta_a\eta_e^{-1}$ in the stalk of an edge with head $a$. The pinned sheaf $\mathcal{G}_A$ acts by left multiplication, trivially on the singletons, and $\torsor_{\eta,\beta}$ is a $\mathcal{G}_A$-torsor; for trivial transport its underlying sheaf of sets is the sheaf of local lifts of $\beta$ along the quotient $\underline{G} \to i_*\underline{G}|_A$. A global section of $\torsor_{\eta,\beta}$ is exactly an interpretation of the unpinned problem, $s_v = s_u\eta_e$ on every edge, that agrees with the cues, $s_a = \beta_a$ on $A$.

Pinned torsors with trivial transport are classified directly. An isomorphism $\torsor_{1,\beta} \to \torsor_{1,\beta'}$ is right multiplication $g \mapsto gh_x$ on each stalk outside $A$ (as in the proof of Theorem \ref{thm:network_torsor_classification}); compatibility with the identity restriction maps makes $h_x = h$ constant on a connected $X$, and matching the pinned singletons forces $\beta'_a = \beta_a h$ for every $a \in A$. Hence $[\torsor_{1,\beta}] = [\torsor_{1,\beta'}]$ precisely when $\beta' = \beta h$ for a single $h \in G$, and $[\torsor_{1,\beta}]$ is the distinguished class precisely when $\beta$ is constant. Write $\partial\beta = [\torsor_{1,\beta}]$ for this boundary-torsor assignment. It is the connecting map of the short exact sequence $1 \to \mathcal{G}_A \to \underline{G} \to i_*\underline{G}|_A \to 1$ in our conventions; for abelian coefficients, under the identification of Proposition \ref{prop:pinned-relative}, it is $\partial = -\delta^*$, the negative of the connecting map of \eqref{eq:LES}. The exact sequence of pointed sets \cite{serre1997galois}
\[
1 \to H^0(X; \mathcal{G}_A) \to H^0(X; \underline{G}) \to H^0(A; \underline{G}) \xrightarrow{\partial} H^1(X; \mathcal{G}_A) \to H^1(X; \underline{G})
\]
holds in the sense that the image of each map is the preimage of the distinguished point under the next. In our left-action convention the global sections $H^0(X;\underline{G}) = G$ act on boundary data by right multiplication, and the fibers of $\partial$ are exactly their orbits.

{\em The Pinned Interval.} Let $X = P_n$ with $A = \{v_0, v_n\}$. A section of $\mathcal{G}_A$ equals $1$ at the pinned ends, and the restriction maps propagate this value along the path, so $H^0(X;\mathcal{G}_A) = \{1\}$. Gauge transformations fixed to $1$ at the ends can move every edge label onto a single edge, and their product along the path is gauge invariant, so $H^1(X,A;\underline{G}) \cong G$ by the ordered product of labels. Under this identification the pinned torsor $\torsor_{\eta,\beta}$ has class $\beta_{v_0}\,\mathrm{hol}_\eta(\gamma)\,\beta_{v_n}^{-1}$, where $\gamma$ is the path from $v_0$ to $v_n$, and it has a global section precisely when
\[
\mathrm{hol}_\eta(\gamma)^{-1}\,\beta_{v_0}^{-1}\,\beta_{v_n} = 1 ,
\]
equivalently $\beta_{v_n} = \beta_{v_0}\,\mathrm{hol}_\eta(\gamma)$. With trivial transport, $H^0(X;\underline{G}) = G$ acts on $G \times G$ by right multiplication on both factors, $H^1(X;\underline{G})$ is trivial, and the classification above gives a bijection
\[
H^1(P_n, \partial P_n; \underline{G}) \;\cong\; (G \times G)/\Delta_G \;\cong\; G, \qquad [(g_0, g_1)] \mapsto g_0\, g_1^{-1} .
\]
This is a bijection of pointed sets: the diagonal $\Delta_G$ need not be normal, so the orbit set is not a quotient group under componentwise multiplication. Nor is it the set of conjugacy classes of $G$. Identifying the two pinned ends turns the interval into a circle, and for abelian $B$ and nonempty $A$ one has $H^1(X,A;B) \cong \widetilde{H}^1(X/A;B)$: a pinned interval is a loop with a marked point. For nonabelian $G$ the loop forgets more, since $H^1(S^1;\underline{G})$ is the set of conjugacy classes of $G$; the pins retain the framing that a loop discards.

{\em Examples.} For the Necker sequence, $G = \mathbb{Z}_2$ and $\beta = (+1,-1)$: the class $g_0 g_1^{-1} = -1$ is the nontrivial element, and no consistent interpretation exists. The impossible bars of Figure \ref{fig:boundary} are the same problem: each end is forced to one of two readings, and the continuity of the bar demands that they agree. A rotational version, with $G = SO(3)$ and distinct end frames $g_0 \neq g_1$, is obstructed only if the bar is modeled as rigid, with identity transport along it; without that constraint a frame may turn continuously along the bar, and nothing is obstructed. In every case the model certifies that no consistent interpretation exists. It does not predict where along the figure the percept switches, nor how ambiguous the middle appears.

%------------------------------------------
\subsection{Pinned Paradoxes on Trees}
%------------------------------------------

With constant coefficients, a tree has no unpinned torsor obstruction; with nonconstant structure sheaves it can, since a single edge whose stalks are all $\mathbb{Z}$, with both restriction maps multiplication by $2$, has $H^1 \cong \mathbb{Z}/2$. Let $X$ be a tree with $N \geq 2$ leaves $L$, pinned at the leaves, with constant coefficients. Boundary assignments that differ by a global section -- in the bistable case, a simultaneous flip of every leaf -- already define the same class in $H^1(X, L; \underline{G})$. We further regard two classes as equivalent when they differ by a combinatorial automorphism of $X$; such an automorphism permutes the leaves and acts on $H^1(X, L; \underline{G})$. We do not identify boundary conditions related only by automorphisms of the coefficient group, which would be a further, separate choice.

\begin{proposition}
\label{prop:bistable-trees}
Let $X$ be a tree with $N \geq 2$ leaves $L$. Then
\[
H^1(X,L;\mathbb{Z}_2) \cong (\mathbb{Z}_2)^N / \langle(-1,\ldots,-1)\rangle \cong (\mathbb{Z}_2)^{N-1},
\]
and the classes of non-trivial bistable paradoxes on $X$, up to automorphisms of $X$, are the orbits of $\Aut(X)$ on the non-trivial elements of this group.
\end{proposition}

{\em Proof:} In the sequence \eqref{eq:LES} for the pair $(X,L)$ with $\mathbb{Z}_2$ coefficients, $H^1(X;\mathbb{Z}_2)$ is trivial because $X$ is a tree, so $\delta^*$ is onto. The group $H^0(L;\mathbb{Z}_2) \cong (\mathbb{Z}_2)^N$ consists of all boundary assignments, and the image of $\res$ consists of the two constant assignments, the subgroup generated by $(-1,\ldots,-1)$. The trivial class is fixed by $\Aut(X)$, so the non-trivial classes are the remaining orbits. \qed

\begin{example}[Star Tree Paradoxes]
For a star $X = T_N$ with $N$ leaves, $\Aut(T_N) \cong S_N$ permutes the leaves arbitrarily. A non-trivial class is a partition of the leaves into two non-empty sets, determined up to exchanging the two sets, and two partitions lie in one orbit precisely when their smaller parts have the same size. The number of classes is therefore $\lfloor N/2 \rfloor$. For $N = 4$ there are two: one leaf against three, and two against two.
\end{example}

%%%%%%%%%%%%%%%%%%%%%%%%%%%%%%%%%%%%%%%%%%%%%
\section{Changing the Model}
\label{sec:isomorphism}
%%%%%%%%%%%%%%%%%%%%%%%%%%%%%%%%%%%%%%%%%%%%%

A model of a figure is a base $X$, a structure sheaf $\mathcal{G}$, and a torsor, recorded by its class in $H^1(X;\mathcal{G})$. The class may be trivial: a comparison of models that admitted only paradoxical ones would be unable to say when a paradox disappears. Three operations relate models. The first is an equivalence; the other two are not, and each can destroy an obstruction. \emph{An obstruction belongs to a model, not to a picture.}

{\em Isomorphism.} For a fixed base and a fixed structure sheaf, torsors are classified up to isomorphism by $H^1(X;\mathcal{G})$ (Theorems \ref{thm:network_torsor_classification} and \ref{thm:nonconstant_classification}). Reindexing by an automorphism of the base or of the coefficients is a further choice, and it should be declared. Reversing the orientation of a loop replaces a holonomy by its inverse; an automorphism $\phi$ of $G$ replaces $[\eta]$ by $[\phi\circ\eta]$. For the cubic staircases of Section \ref{sec:cubic_staircases}, each holonomy is $k$ times a primitive vector, so under the full group $GL(3,\mathbb{Z})$ of lattice automorphisms the eight examples fall into three classes, according as $|k|$ is $1$, $2$, or $4$. Arbitrary lattice automorphisms do not preserve the geometry of unit cubes, however. Under the signed permutations of the axes, the symmetries of the cubical lattice, $(1,1,1)$ and $(2,2,1)$ are inequivalent, and the eight examples fall into four classes: $(1,1,1)$, $(2,2,2)$, $(4,4,4)$, and $(2,2,1)$. The two rows of Figure \ref{fig:penrose-stairs} are different paradoxes.

{\em Pullback.} A map $f: Y \to X$ pulls torsors back,
\[
f^* : H^1(X;\mathcal{G}) \to H^1(Y; f^{-1}\mathcal{G}),
\]
and for constant coefficients $f^*$ precomposes the holonomy with $f_*: \pi_1(Y) \to \pi_1(X)$. The inclusion of the first loop of the torus graph pulls the torus staircase with periods $(1,0)$ back to the Penrose staircase with climb $1$; the inclusion of the second loop pulls the same class back to the trivial one. Pullback along the two loops of the Klein ladder exhibits its parts: along $b$ a translation $(\beta,+1)$, a cylindrical staircase; along $a$ a reflection $(\alpha,-1)$, a M\"obius staircase. For the pinned sheaves of Section \ref{sec:boundary}, pullback exists along any map but transports the pinning to $f^{-1}(A)$: for maps compatible with the cellular structures, $f^{-1}\mathcal{G}_A$ is the sheaf pinned on $f^{-1}(A)$, so identifying it with a specified $\mathcal{G}_B$ requires $B = f^{-1}(A)$. A map of pairs, which asks only $B \subseteq f^{-1}(A)$, is not enough for that identification, though it still induces a map of relative cohomology. Nor is the pinned sheaf homotopy invariant: two maps from a point into a path, homotopic through the path, pull it back to the trivial group at an endpoint and to $G$ in the interior.

{\em Change of Coefficients.} A homomorphism of structure sheaves $\phi: \mathcal{G} \to \mathcal{H}$ pushes classes forward,
\[
\phi_* : H^1(X;\mathcal{G}) \to H^1(X;\mathcal{H}),
\]
by applying $\phi$ to a cocycle. Some coefficient changes lose nothing that matters. On a circle, $\Phi: \mathbb{Z}^3 \to \mathbb{Z}$, $\Phi(a,b,c) = -a + b$, carries the cubic holonomy $(2,4,6)$ to the height holonomy $2$, and $\Psi: \mathbb{Z} \to \mathbb{Z}^3$, $\Psi(n) = (n, 2n, 3n)$, carries it back, although neither map is invertible. Others forget. The sign homomorphism $D_\infty \to \mathbb{Z}_2$ carries the M\"obius staircase to the zigzag paradox, keeping the flip and discarding the height. Reduction modulo two carries a Penrose staircase with climb $k$ to the $\mathbb{Z}_2$-class $(-1)^k$, and so kills every staircase with even climb.

Pullback and change of coefficients preserve the distinguished class, and both can carry a non-trivial class to it. Neither defines an equivalence of reconstruction problems. Relating two figures by such maps says which information one model retains about the other; it does not say that the figures are the same paradox.

%%%%%%%%%%%%%%%%%%%%%%%%%%%%%%%%%%%%%%%%%%%%%
\section{Closing the Loop on the Penrose Triangle}
\label{sec:closeloop}
%%%%%%%%%%%%%%%%%%%%%%%%%%%%%%%%%%%%%%%%%%%%%

The Penrose triangle of Figure \ref{fig:classics} invites several models, and they are not interchangeable. We give one in full, with explicit states and transport, and then place it in a setting that also contains the M\"obius staircase and the Klein ladder.

{\em The Tribar as a Rigid Frame.} Model the tribar by a cycle with three vertices, the corners, and three edges, the beams. A local state at a corner is a rigid placement of that corner in space: a position together with an orientation. The group $SE(3)$ of rigid motions, with product $(R,\vect{t})(S,\vect{u}) = (RS, \vect{t} + R\vect{u})$, acts simply transitively on such placements, so the local states form an $SE(3)$-torsor. Each corner is read as a standard right-angle joint in a common set of axes, and each beam is straight. The transition along the $i$-th beam is therefore a pure translation $(I, L\vect{v}_i)$, with $\vect{v}_i$ a unit axis vector and $L$ the length of a beam, and the holonomy is
\[
\bigl(I,\; L(\vect{v}_1 + \vect{v}_2 + \vect{v}_3)\bigr).
\]
The beams run along the three axes and the drawing closes, so, as for the cubic staircases of Section \ref{sec:cubic_staircases}, the translation lies in the kernel of the isometric projection: in adapted coordinates it is $L(1,1,1) \neq 0$. The rotational part of the holonomy is the identity. The three right-angle turns are turns of the beam's tangent, not rotations of the frame; carried around the triangle, the frame returns unrotated and displaced. Since the holonomy lies in the abelian subgroup of translations, the torsor reduces to an $\mathbb{R}^3$-torsor, and projection to the vertical coordinate, a homomorphism on translations, gives a height class $\pm L \neq 0$. Under this reduction to translations, the tribar carries the scalar height obstruction of a Penrose staircase. \emph{The triangle hides a staircase.}

{\em Other Models.} A rotational model with structure group $SO(3)$ would need transitions that rotate the frame, and none is visible in the figure: the projection $SE(3) \to SO(3)$ carries the class above to the trivial one. Penrose's multiplicative model \cite{Penrose1992Cohomology}, with coefficients in $\mathbb{R}^+$, is identified with the height model by the isomorphism $\exp: (\mathbb{R},+) \to (\mathbb{R}^+,\times)$, when the holonomies correspond. A binary model of front and behind would need its own local states and transitions: a cyclic ordering of occlusions is not by itself a $\mathbb{Z}_2$-torsor, and we do not pursue it here.

These models cannot be traded for one another freely. The groups $SO(3)$ and $SE(3)$ are perfect. Every rotation is a commutator, $[R_{\theta/2}, S] = R_\theta$ for a half-turn $S$ about an axis perpendicular to that of $R_\theta$, and every translation is a commutator, $[(S,\vect{0}),(I,\vect{v})] = (I, (S-I)\vect{v})$, which equals $(I,\vect{t})$ when $S$ is a half-turn reversing $\vect{t}$ and $\vect{v} = -\vect{t}/2$. Neither group therefore admits a nontrivial homomorphism to an abelian group, and no coefficient map carries a nontrivial rigid-motion class to a nontrivial height, scale, or binary class. In particular, extracting the translational part $(R,\vect{t}) \mapsto \vect{t}$, or equivalently evaluating at the origin, is not a homomorphism. The height class of the tribar arises by reduction to the translation subgroup, available because this particular holonomy is a translation, and not by a coefficient map out of $SE(3)$.

{\em Affine Holonomy.} When the frames do rotate, the translation part survives in a twisted form. Let $V$ be an abelian group or vector space on which a group $H$ acts linearly, and let
\[
\rho : \pi_1(X) \to V \rtimes H, \qquad \rho(\gamma) = (R(\gamma), \vect{t}(\gamma)).
\]
Multiplying out gives $\vect{t}(\gamma\delta) = \vect{t}(\gamma) + R(\gamma)\,\vect{t}(\delta)$: the translation part is a crossed homomorphism for the local coefficient system $V_R$ determined by the linear part. Moving the affine origin to $\vect{c}$ conjugates by $(I,-\vect{c})$ and replaces $\vect{t}(\gamma)$ by $\vect{t}(\gamma) + (R(\gamma) - I)\vect{c}$, a twisted coboundary. For a fixed linear part, the translations therefore define a class
\[
[\vect{t}] \in H^1(X; V_R),
\]
the obstruction to a compatible affine origin. It is not the whole frame torsor: $[\vect{t}] = 0$ says that the origins can be made consistent once the linear transport is fixed, not that the linear transport is trivial. For a single loop,
\[
H^1(S^1; V_R) \cong V/(R - I)V .
\]
Three cases recover the examples of this paper. When $R = I$, as for the tribar, the cubic staircases, and the Penrose staircase, the quotient is $V$ itself and the translation survives intact. When $R$ is a nonidentity rotation of $\mathbb{R}^3$, the quotient is the axis of $R$: the axial displacement of a screw motion is exactly what no change of origin removes, which is why projecting a rigid-motion class to $SO(3)$ and including it back loses the class precisely when the axial displacement is nonzero. When $V = \mathbb{Z}$ and $R = -1$, the quotient is $\mathbb{Z}/2$: this is the parity $\Delta \bmod 2$ of the M\"obius staircase in Section \ref{sec:non_orientable}, and the conjugation formula there is the twisted coboundary here.

The Klein ladder is the two-loop case. Its linear part is the orientation character $w$ of the Klein bottle, with $w(a) = -1$ and $w(b) = +1$. A crossed homomorphism on $F_2 = \pi_1(\Lambda)$ is an arbitrary pair $(\vect{t}(a), \vect{t}(b)) = (\alpha, \beta)$; every such pair respects the relator $aba^{-1}b$, as computed in Section \ref{sec:nonabelian_paradox}; and the twisted coboundaries are the pairs $(-2c, 0)$. Hence
\[
H^1(K; \mathbb{Z}_w) \cong \mathbb{Z}/2 \oplus \mathbb{Z}, \qquad [\vect{t}] \mapsto (\alpha \bmod 2, \beta),
\]
which by twisted Poincar\'e duality is $H_1(K;\mathbb{Z})$. The invariants of the Klein ladder found in Section \ref{sec:nonabelian_paradox} are this class, up to the sign of $\beta$, which a reflection of the height line reverses. The height paradoxes of this paper, straight and twisted, are affine holonomy with coefficients in a line.

%%%%%%%%%%%%%%%%%%%%%%%%%%%%%%%%%%%%%%%%%%%%%
\section{Conclusion}
\label{sec:conclusion}
%%%%%%%%%%%%%%%%%%%%%%%%%%%%%%%%%%%%%%%%%%%%%

Penrose's seminal observation linking the impossible triangle to the first cohomology group of the circle points to a deep connection between visual paradoxes and obstruction theory. This paper contributes a precise mathematical framework for analyzing paradoxical figures.

At the heart of this framework lies the network torsor -- a discrete combinatorial analogue of a classical torsor that captures the structure of visual paradox. Once a model specifies its base graph, its structure group, its local states, and the transport between them, the impossibility of a figure is precisely the obstruction to finding a global section. When the local states form a torsor, this obstruction is quantified by a non-trivial element of the first cohomology $H^1(X; \mathcal{G})$. The framework accommodates both loop-based paradoxes (classified by cohomology with constant coefficient sheaves) and boundary-driven paradoxes (classified by pinned structure sheaves, equivalently by relative cohomology), providing a unified perspective across diverse paradoxical phenomena.

While this paper has focused on paradoxes classified by first cohomology over graphs, surfaces enter only through conditions on the holonomy around their 2-cells. Obstructions carried by the 2-cells themselves, in second cohomology, are considered briefly in \cite{GhristGhrist2026Degrees}. Also of potential interest is the degree to which non-constant structure sheaves arise in visual paradoxes beyond boundary conditions and the twisted coefficients of Section \ref{sec:closeloop}.

Additionally, our presentation of network torsors connects to current research in applied topology and distributed systems. The exploration of the relationship between torsors and other sheaf-theoretic constructions over discrete spaces may yield both theoretical insights and practical algorithms for detecting and analyzing inconsistencies in geometric data.

%%%%%%%%%%%%%%%%%%%%%%%%%%%%%%%%%%%%%%%%%%%%%
%%%%%%%%%%%%%%%%%%%%%%%%%%%%%%%%%%%%%%%%%%%%%
%%%%%%%%%%%%%%%%%%%%%%%%%%%%%%%%%%%%%%%%%%%%%

\section*{Acknowledgments}
This work was initially funded by the Air Force Office of Scientific Research under award number FA9550-21-1-0334. Research and writing was done in collaboration with AI language models; all artwork is original and without use of AI.  
Certain images of generalized paradoxes were first created in March 2022, such as Figure \ref{fig:torus}[right]. 
%The writing of this paper was begun on March 23, 2025. 

\bibliographystyle{plain}
\bibliography{IMPOSSIBLE}

\end{document}